\documentclass[a4paper,11pt]{article}
\usepackage{amsfonts,amsmath,amssymb,amsthm,verbatim,placeins,comment,placeins,caption,subcaption,enumerate,bm,graphicx,fancyvrb}
\usepackage[dvipsnames]{xcolor}
\usepackage{geometry}
\usepackage{float}
\usepackage[pdfauthor={DG},pdfstartview={FitH},colorlinks=true,citecolor=blue,urlcolor=blue,linkcolor=blue]{hyperref}
\usepackage[english]{babel}
\usepackage[sort&compress]{natbib}

\usepackage{setspace}
\newtheorem{theorem}{Theorem}

\newtheorem{proposition}{Proposition}
\newtheorem{definition}{Definition}
\newtheorem{example}{Example}
\theoremstyle{remark}

\newcommand{\E}{\mathbb{E}}

\newcommand{\R}{\mathbb{R}}
\newcommand\1{\mathbf{1}}
\newcommand{\I}{\operatorname{I}}
\newcommand{\II}{\operatorname{II}}
\newcommand{\III}{\operatorname{III}}
\newcommand\D{\mathcal{D}}
\newcommand\A{\mathcal{A}}

\DeclareMathOperator{\Leb}{Leb}
\newcommand{\Var}{\mathop{\mathrm{Var}}}
\newcommand{\Cov}{\mathop{\mathrm{Cov}}}
\newcommand{\Beta}{\operatorname{Beta}}
\newcommand{\B}{\operatorname{B}}
\renewcommand{\Im}{\operatorname{Im}}
\renewcommand{\Re}{\operatorname{Re}}

\usepackage{booktabs}
\usepackage{bera}
\usepackage{xurl}
\usepackage{appendix}

\usepackage{authblk}

\makeatletter
\let\@fnsymbol\@alph
\makeatother

\title{Superpositions of CARMA processes}
\date{}
\author[1]{Danijel Grahovac\thanks{dgrahova@mathos.hr}}
\author[1]{Magdalena Mikić\thanks{mmikic1@mathos.hr}}

\affil[1]{School of Applied Mathematics and Informatics, J.~J.~Strossmayer University of Osijek, Trg Ljudevita Gaja 6, 31000 Osijek, Croatia}

\begin{document}

\maketitle

\smallskip

\noindent \textbf{Abstract: } We introduce supCARMA processes, defined as superpositions of L\'evy-driven CARMA processes with respect to a L\'evy basis, as a natural extension of the superpositions of Ornstein-Uhlenbeck type processes. We then focus on supCAR$(2)$ processes and show that they can be classified into three distinct types determined by the eigenstructure of the underlying CAR$(2)$ matrix. For each type we provide conditions for existence and derive explicit expressions for the correlation function. The resulting correlation structures may exhibit long-range dependence and can be non-monotone. These features make supCAR$(2)$ processes a flexible class for modeling time series with oscillatory correlations or strong dependence.

\bigskip


\bigskip

\noindent\textbf{Funding: } This work was supported by the Croatian Science Foundation (HRZZ) grant Scaling in Stochastic Models (IP-2022-10-8081).


\noindent\textbf{Disclosure statement: } The authors report there are no competing interests to declare.

\bigskip

\section{Introduction}

The idea that aggregation of simple dynamic models may generate rich dependence structures goes back to the seminal work of Granger (see \cite{granger1980}), who showed that contemporaneous aggregation of short-memory autoregressive processes of order 1 (AR$(1)$) can produce time series with long-range dependence and polynomially decaying (auto)corre\-la\-tions. Here, long-range dependence means that the correlation function is not integrable. Granger considered AR$(1)$ processes with coefficients drawn from a distribution that places substantial mass near one, together with independent finite-variance noise, and demonstrated that their aggregate can display much slower decay of correlations than the individual components. Related ideas also appear earlier in \cite{robinson1978statistical}, who analyzed random-coefficient AR$(1)$ models not to study the aggregation itself, but to identify the distribution of the autoregressive coefficient from cross-sectional autocovariances. 

Following Granger's contribution, numerous studies have shown that long-range dependence may arise through the aggregation of random-coefficient autoregressive models. In \cite{linden1999}, the aggregation of AR$(1)$ processes with uniformly distributed coefficients was considered. Extensions to AR$(2)$ models were considered in \cite{chong2001}, where the aggregation of AR$(2)$ processes with randomly distributed characteristic roots was considered. Aggregation of higher order autoregressive processes has been considered in \cite{oppenheim2004,Zaffaroni2004}. Although the aggregated processes in all these works are Gaussian, a different type of aggregation scheme was introduced in \cite{philippe2014}, who considered a triangular array of random-coefficient AR$(1)$ processes with innovations in the domain of attraction of an infinitely divisible law. In this setting, the aggregate need not be Gaussian. The limit distribution depends on both the innovations and how the coefficients behave near one. It can be Gaussian, but it may also be non-Gaussian and infinitely divisible, including stable laws. This extends the results of \cite{puplinskaite2009, puplinskaite2010}, in which stable limits were obtained in the infinite-variance case.

Earlier work in continuous time was carried out by Robinson (see \cite{robinson1986cross}), who analyzed continuous-time analogs of AR models of first and second order with randomly chosen autoregressive coefficients and finite-variance noise. He derived explicit expressions for the autocovariance function of the aggregated process and showed that its dependence is governed by the distribution of these random coefficients.
A further development in continuous time was introduced by Barndorff-Nielsen (see \cite{bn2001}), who proposed superpositions of Ornstein-Uhlenbeck (OU) type processes driven by a L\'evy basis, known as supOU processes. The use of L\'evy noise allows jumps, heavy tails, and infinitely divisible marginals. Moreover, the marginal distribution and the correlation function can be modeled independently, and it is possible to produce long-memory models (see, e.g.~\cite{bnleonenko2005}). Moreover, supOU processes are infinitely divisible and fall into the class of ambit processes (see \cite{barndorff2018ambit}).  The construction of supOU processes has been extended to the multivariate setting in \cite{barndorff2011multivariate}. SupOU processes also exhibit interesting limiting behavior (see \cite{GLST2019,GLT19,GLT2020multifaceted,grahovackevei2025}). Recently, a spatial extension has been developed in \cite{donhauzer2025construction}, where supCAR fields are introduced. They are defined as superpositions of CAR$(1)$ random fields with respect to a L\'evy basis. 

The OU type processes are L\'evy-driven continuous time autoregressive processes of order 1 (CAR$(1)$)  \cite{brockwell2001249,brockwell2001levy}. Hence, the supOU processes can be considered as supCAR$(1)$ processes. The aim of this paper is to generalize the supOU framework to the continuous time ARMA setting by introducing supCARMA processes. We define these processes as integrals with respect to infinitely divisible random measures and establish conditions under which these integrals are well-defined. This yields a flexible class of mixed moving average processes in which the dependence structure is determined by the distribution of the CARMA parameter matrices. We focus on the case of supCAR$(2)$ processes. Based on the eigenstructure of the underlying CAR$(2)$ matrix, we show that three structurally different types of supCAR$(2)$ processes arise. For each type, we establish existence and derive explicit expressions for the correlation function. In addition to the possibility of long-range dependence, the correlation function may also exhibit oscillatory behavior. This contrasts the supOU processes whose correlation function is always positive and decreasing.

The superpositions of CARMA processes have already appeared in unpublished manuscript \cite{marquardt2007generating}. Unlike the scalar randomization used in \cite{marquardt2007generating}, our construction of supCARMA processes follows the perspective of Barndorff-Nielsen and Stelzer in \cite{barndorff2011multivariate}, where supCAR models appear as linear functionals of multivariate supOU processes. 

SupOU processes have been applied in a wide range of settings that require flexible dependence structures and infinitely divisible noise. They have proven particularly useful in financial econometrics and stochastic volatility modeling (see \cite{barndorff2013multivariate, barndorff2001nongaussian, moser2011tail, stelzer2015moment}), and also in the modeling of turbulence and other physical or environmental systems through their connection to ambit processes (see \cite{barndorff2018ambit}). Their ability to combine prescribed marginal distributions with flexible correlation structures makes supOU processes particularly suitable for modeling time series with jumps, heavy tails, or long-range dependence. Within this context, the supCARMA framework, and in particular the supCAR$(2)$ class introduced in this paper, offers additional modeling potential by allowing mixtures of second-order dynamics. This enables the representation of oscillatory or non-monotonic correlation functions and makes supCAR$(2)$ processes promising candidates for modeling complex time series with oscillatory dependence or long memory behavior.

In Section \ref{sec:preliminaries} we recall basic definitions and notation, including L\'evy bases and L\'evy-driven CARMA processes. Section \ref{sec:supCARMA} introduces supCARMA processes and establishes their representation both as mixed moving averages and as linear functionals of multivariate supOU processes. In Section \ref{sec:supCAR2} we specialize to the case $p=2$ and obtain the decomposition of supCAR$(2)$ processes into three types based on the discriminant of the characteristic polynomial. Sections \ref{sec:typeI}, \ref{sec:typeII} and \ref{sec:typeIII} analyze the three types separately. For each type, we provide integrability conditions ensuring existence, derive moment conditions, and study second-order properties including the correlation function.

\section{Preliminaries}\label{sec:preliminaries}

\subsection{Notation}

In the following, $|\cdot|$ denotes the Euclidean norm of the vector and $\|\cdot \|$ the matrix norm, which, if not specified, we take to be the spectral norm. $C$ denotes the constant that may change from row to row.

We will denote the cumulant (generating) function of a (multivariate) random variable $Y$ by
\begin{equation*}
\kappa_Y(\bm{z})= \log \E e^{i \bm{z}^T Y}.
\end{equation*}
For a stochastic process $Y=\{Y(t)\}$, we write $\kappa_Y(z,t) = \kappa_{Y(t)}(z)$, and for the cumulant function of the random variable $Y(1)$ we write $\kappa_Y(z)=\kappa_Y(z,1)$. 

If $Y$ is infinitely divisible, then by the L\'{e}vy-Khintchine formula 
\begin{equation*}
\kappa_Y(\bm{z}) =i\gamma^T \bm{z} -\frac{1}{2}\bm{z}^T \Sigma \bm{z} + \int_{\R^d}\left( e^{i\bm{z}^T \bm{x}}-1-i\bm{z}^T \bm{x} \mathbf{1}_{[0,1]}(|\bm{x}|)\right) \mu(d\bm{x}), \quad \bm{z} \in \R^d,
\end{equation*}
where the characteristic triplet $(\gamma,\Sigma,\mu)$ consists of $\gamma\in \R^d$, $\Sigma$ a symmetric nonnegative definite $d\times d$ matrix and L\'evy measure $\mu$ on $\R^d$ such that $\mu\left( \left\{ 0\right\} \right) =0$ and $\int_{\R^d} (1 \wedge |\bm{x}|^2) \mu(d\bm{x})<\infty$.
For a measure $\pi$ on $(0,\infty)$ and $p \in \R$, we denote by
\begin{equation*}
    m_p(\pi) = \int_0^{\infty} a^p \, \pi(da)
\end{equation*}
the $p$-th moment of $\pi$ whenever the integral exists.

\subsection{CARMA processes}
A zero-mean L\'evy-driven CARMA$(p,q)$ process $Y$ with $0\leq q<p$ can be defined by the state-space representation
\begin{equation}\label{e:CARMA}
Y(t) = \bm{b}^T \bm{X}(t),
\end{equation}
where $\bm{X}(t)$ is a solution of the stochastic differential equation (SDE)
\begin{equation}\label{e:SDEcarmaX}
d \bm{X}(t) = A \bm{X}(t) dt + \bm{e} dL(t)
\end{equation}
with
\begin{equation}\label{e:Abe}
A= \begin{bmatrix}
0 & 1 & 0 & \cdots & 0\\
0 & 0 & 1 & \cdots & 0\\
\vdots & \vdots & \vdots & \ddots & \vdots\\
0 & 0 & 0 & \cdots & 1\\
-a_p & -a_{p-1} & -a_{p-2} & \cdots & -a_{1}
\end{bmatrix}, \quad
\bm{e}=\begin{bmatrix}
0\\
0\\
\vdots\\
0\\
1
\end{bmatrix}, \quad
\bm{b}= \begin{bmatrix}
b_0\\
b_1\\
\vdots\\
b_{p-2}\\
b_{p-1}
\end{bmatrix},
\end{equation}
$b_j=0$ for $q<j\leq p$ and $L$ being a L\'evy process. A solution of the SDE \eqref{e:SDEcarmaX}, given $\bm{X}(0)$ independent of $L$, is 
\begin{equation*}
\bm{X}(t) = e^{At} \bm{X}(0) + \int_0^t e^{A(t-s)} \bm{e} dL(s).
\end{equation*}
Formally, one can view CARMA process as a stationary solution of $p$-th order linear differential equation
\begin{equation*}
a(D) Y(t) = b(D) DL(t), \quad t \geq 0,
\end{equation*}
where $L$ is some L\'evy process, $D$ denotes differentiation with respect to $t$, and
\begin{align*}
a(z) &= z^p + a_1 z^{p-1} + \cdots +a_{p-1}z+ a_p,\\
b(z) &= b_0 + b_1 z + \cdots + b_{p-1}z^{p-1} + b_p z^p.
\end{align*}

\medskip

If $\E \log (1+|L(1)|)<\infty$ and the eigenvalues of $A$ (or, equivalently, the roots of the polynomial $a$) all have strictly negative real parts, then SDE \eqref{e:SDEcarmaX} has a strictly stationary solution given by
\begin{equation}\label{e:OUdef}
\bm{X}(t) = \int_{-\infty}^t e^{A(t-s)} \bm{e} dL(s).
\end{equation}
Here, $L$ is extended to a L\'evy process on the whole real line by letting $\{\widetilde{L}(t), \, t\geq 0\}$ be an independent copy of $L$, and defining $L(t)=-\widetilde{L}(-t-)$ for $t<0$. Moreover, CARMA$(p,q)$ process $Y$ is a unique strictly stationary solution satisfying \eqref{e:CARMA} and \eqref{e:OUdef} (\cite{brockwelllindner2009}; see also \cite{brockwell2005levy,brockwell2001249,brockwell2001levy}). We can write $Y$ in the moving average representation as
\begin{equation*}
Y(t) = \bm{b}^T \int_{-\infty}^t e^{A(t-s)} \bm{e} dL(s) = \int_{-\infty}^t \bm{b}^T e^{A(t-s)} \bm{e} dL(s).
\end{equation*}

\subsection{L\'evy bases}

Let $(S, \mathcal{S})$ be a Borel (Lusin) space with $\mathcal{S}$ a Borel $\sigma$-algebra on $S$. We denote by $\mathcal{B}_b(S)$ the bounded sets in $\mathcal{S}$. An $\R^d$-valued L\'evy basis $\Lambda$ on $(S, \mathcal{S})$ is an infinitely divisible independently scattered random measure, i.e.~a collection of random variables $\{\Lambda(A), \, A \in \mathcal{B}_b(S)\}$ such that
\begin{enumerate}[(i)]
\item For any sequence $A_1,A_2,\dots$ of disjoint elements of $\mathcal{B}_b(S)$ such that $\cup_{i=1}^\infty A_i \in \mathcal{B}_b(S)$ we have $\Lambda \left( \cup_{i=1}^\infty A_i\right) = \sum_{i=1}^\infty \Lambda(A_i)$ a.s.
\item For any sequence $A_1,A_2,\dots$ of disjoint elements of $\mathcal{B}_b(S)$, $\Lambda(A_1), \Lambda(A_2),\dots$ are independent.
\item For any $A \in \mathcal{B}_b(S)$, the distribution of $\Lambda(A)$ is infinitely divisible.
\end{enumerate}

We consider only homogeneous L\'evy bases and in this case, for any $A \in \mathcal{B}_b(S)$,
\begin{equation*}
\E e^{i \bm{z}^T \Lambda(A)} = e^{\kappa(\bm{z}) M(A)}, \quad \bm{z} \in \R^d,
\end{equation*}
where $M$ is a measure on $(S, \mathcal{S})$ (control measure) and $\kappa$ is the characteristic exponent of some infinitely divisible distribution determined by some triplet $(\gamma,\Sigma,\mu)$. The distribution of $\Lambda$ is completely determined by the control measure $M$ and L\'evy-Khintchine triplet $(\gamma,\Sigma,\mu)$. See \cite{barndorff2018ambit,samorodnitsky2016book} for more details on L\'evy bases.

The integration of a function on $S$ with respect to the random measure $\Lambda$ can be defined first for real simple functions, then as a limit in probability of such integrals. The conditions for integrability of functions with respect to $\Lambda$ can be found in \cite{rajput1989spectral}; see also \cite{barndorff2018ambit} and the references therein. In particular, for the homogeneous L\'evy basis $\Lambda$, a measurable function $f \colon S \to \R$ is integrable with respect to $\Lambda$ if and only if the following conditions hold
\begin{equation}\label{eq:rajros:gen}    
\begin{aligned}
& \int_S \left| \gamma f(s) + \int_{\R} \left( \tau(x f(s)) - f(s) \tau(x) \right) \mu(dx) \right| M(ds) <\infty,\\
& \int_S \left| \Sigma f(s) \right|^2 M(ds) < \infty,\\
& \int_S \left( \int_{\R} \left( 1 \wedge |x f(s)|^2 \right) \mu(dx) \right) M(ds) < \infty,
\end{aligned}
\end{equation}
where $\tau(x) = x \1_{[0,1]}(| x |)$.

In the following, the typical space $S$ will be of the form $S=V\times \R$, where $V$ is some Borel space. Moreover, the control measure $M$ will be a product of some (probability) measure $\pi$ on $V$ and the Lebesgue measure on $\R$. In this case, we will say that $(\gamma,\Sigma,\mu,\pi)$ is the generating quadruple of the L\'evy basis $\Lambda$ on $V\times \R$.

\section{SupCARMA processes}\label{sec:supCARMA}

SupCARMA processes can be seen as a generalization of the supOU processes defined in \cite{bn2001}. A supOU process is a stationary process given by
\begin{equation*}
X(t) = \int_{(0,\infty)} \int_{-\infty}^{t} e^{-a(t-s)} \Lambda(da, ds), \quad t\in \R,
\end{equation*}
where $\Lambda$ is an $\R$-valued L\'evy basis on $(0,\infty)\times \R$ with generating quadruple $(\gamma,\Sigma,\mu,\pi)$ for some probability measure $\pi$ on $(0,\infty)$. Since L\'evy-driven OU process corresponds to CAR$(1)$ (CARMA$(1,0)$) process, supCARMA may be defined as follows. 

First, we introduce some notation. Let $\sigma(A)$ denote the spectrum of matrix $A$. By $M_p^-:= \{A \in M_p(\R): \ \sigma(A) \subset (-\infty,0)+i\R \}$ we denote the set of all $p$-dimensional matrices whose eigenvalues have a strictly negative real part and $\A_p$ will denote the set of all matrices in $M_p^-$ of the form \eqref{e:Abe}.

\begin{definition}\label{def:supCARMA}
Let $0\leq q <p$ and let $\Lambda$ be an $\R$-valued L\'evy basis on $\A_p \times \R$ with generating quadruple $(\gamma,\Sigma,\mu,\pi)$. A supCARMA$(p,q)$ process $X$ is defined by
\begin{equation}\label{e:def:supCARMA}
X(t) = \int_{\A_p} \int_{-\infty}^t \bm{b}^T e^{A(t-s)} \bm{e} \Lambda(dA, ds), \quad t \in \R,
\end{equation}
where $\bm{b}$ and $\bm{e}$ are given by \eqref{e:Abe}. 
\end{definition}

Note that supCARMA process is a so-called (causal) mixed moving average process since it can be represented as
\begin{equation*}
X(t) = \int_{V} \int_{\R} f(v, t-s)\, \Lambda(dv,ds),
\end{equation*}
for $V=\A_p$ and
\begin{equation}\label{e:mmaf}
f(A,u) = \bm{b}^T e^{Au} \bm{e} \1_{(0,\infty)}(u).
\end{equation}

SupCARMA processes were defined differently in \cite{marquardt2007generating}, where the matrix in the kernel function is taken to be of the form $\lambda A$ for fixed matrix $A$ and randomized parameter $\lambda$.

A definition of a multivariate supCAR process that corresponds to Definition \ref{def:supCARMA} has been noted in \cite[Subsection 5.3]{barndorff2011multivariate}. These processes can be regarded as components of multivariate supOU processes. We detail this argument for supCARMA processes.

Let $\widetilde{\Lambda}(dA, ds) = \bm{e} \Lambda(dA, ds)$, which is an $\R^p$-valued L\'evy basis on $\A_p \times \R$ with generating quadruple $(\gamma,\Sigma,\mu,\pi)$ and such that
\begin{equation*}
\E e^{i \bm{z}^T \widetilde{\Lambda}(B)} = e^{\widetilde{\kappa}(\bm{z}) (\pi \times \Leb)(B)}, \quad \bm{z} \in \R^p. 
\end{equation*}
The characteristic exponent $\widetilde{\kappa}$ corresponds to an infinitely divisible distribution with first $p-1$ components degenerate so that for $\bm{z}=(z_1,\dots,z_p)^T$
\begin{equation*}
\widetilde{\kappa}(\bm{z}) =i \gamma z_p -\frac{1}{2} \Sigma z_p^2 + \int_{\R}\left( e^{iz_p x}-1-iz_p x \mathbf{1}_{[0,1]}(|x|)\right) \mu(dx).
\end{equation*}
We can then write \eqref{e:def:supCARMA} in the state-space representation as
\begin{equation*}
X(t) = \bm{b}^T \bm{X}(t),
\end{equation*}
where
\begin{equation*}
\bm{X}(t) = \int_{\A_p} \int_{-\infty}^t e^{A(t-s)} \widetilde{\Lambda}(dA, ds).
\end{equation*}
The process $\bm{X}$ is a $p$-dimensional supOU process as defined in \cite{barndorff2011multivariate}. Hence, any supCARMA process can be represented as a linear functional of a specific multivariate supOU process.

\subsection{Existence from the state-space representation}

For the supCARMA process to be well-defined it is necessary and sufficient that the supOU process $\bm{X}$ is well-defined, which is then clearly stationary. Hence, $X$ is stationary if and only if $\bm{X}$ is stationary. From \cite[Theorem 3.1]{barndorff2011multivariate} we immediately get the sufficient conditions for the existence of a stationary supCARMA process. The necessary conditions for the existence of the multivariate supOU process are given in \cite[Theorem 3.1]{barndorff2011multivariate}, from which we immediately get the following.

\begin{theorem}
Suppose that $\int_{|x|>1} \log |x| \mu(dx)<\infty$ and that there exist measurable functions $\rho \colon \mathcal{A}_p \to (0,\infty)$ and $\eta \colon \mathcal{A}_p \to [1,\infty]$ such that
\begin{equation*}
\|e^{A s}\| \leq \eta(A) e^{-\rho(A) s}, \quad \forall s\geq 0, \ \pi \text {-a.s.}, 
\end{equation*}
and
\begin{equation}\label{e:mvsupOUcond}
\int_{\mathcal{A}_p} \frac{\eta(A)^2}{\rho(A)} \pi(d A)<\infty.
\end{equation}
Then the supCARMA process \eqref{e:def:supCARMA} is well-defined and stationary.
\end{theorem}

The conditions of the previous theorem seem intricate to check in general. In particular cases, bounds for the norm of matrix exponential could be used. For example, by \cite{kaagstrom1977bounds} we have that
\begin{equation*}
\|e^{A s}\|  \leq \| P \| \| P^{-1}\| e^{s (-1-\max \Re(\sigma(A))}, \quad s \geq 0.
\end{equation*}
Assume that $\max \Re(\sigma(A)) <-1$, $\pi$-a.s. If we put $\rho(A)=-\max \Re(\sigma(A))-1$ and $\eta(A)=\|P\|_2 \|P^{-1}\|_2$, then \eqref{e:mvsupOUcond} reads
\begin{equation*}
\int_{\mathcal{A}_p} \frac{\|P\|_2 \|P^{-1}\|_2}{-\max \Re(\sigma(A))-1} \pi(d A)<\infty,
\end{equation*}
which in some cases may be easier to check.


\subsection{Existence from the mixed moving average representation}

For specific cases, it may be more convenient to express the process in the mixed moving average form \eqref{e:def:supCARMA}. In this case, the Rajput-Rosinski conditions \eqref{eq:rajros:gen} can be used to check whether the process is well-defined. For an $\R$-valued L\'evy basis $\Lambda$ on $V \times \R$ with generating quadruple $(\gamma,\Sigma,\mu,\pi)$ and a measurable function $f \colon V \times \R \to \R$, the mixed moving average process
\begin{equation*}
    X(t) = \int_{V} \int_{\R} f(v, t-s)\, \Lambda(dv,ds)
\end{equation*}
is well-defined if and only if
\begin{equation}\label{eq:rajros:mma}
\begin{aligned}
& \int_{V\times \R} \left| \gamma f(v, s) + \int_{\R} \left( \tau(x f(v, s)) - f(v, s) \tau(x) \right) \mu(dx) \right| \pi(dv) ds <\infty,\\
& \int_{V\times \R} \left| \Sigma f(v, s) \right|^2  \pi(dv) ds < \infty,\\
& \int_{V\times \R} \left( \int_{\R} \left( 1 \wedge |x f(v, s)|^2 \right) \mu(dx) \right) \pi(dv) ds < \infty.
\end{aligned}
\end{equation}
For specific examples of the supCARMA process, these conditions should be checked for $f$ given in \eqref{e:mmaf}.

\section{SupCAR$(2)$ processes}\label{sec:supCAR2}

The properties of supCARMA processes depend crucially on the properties of the measure $\pi$. In general, it is far from obvious how to define distribution on the space of matrices $\mathcal{A}_p$ to ensure, e.g., that the process is well-defined or that it has a certain correlation structure. Therefore, we provide a deeper analysis of the $p=2$ and $q=0$ case.

Consider a process obtained from Definition \ref{def:supCARMA} by taking $p=2$ and $q=0$
\begin{equation*}
		X(t) = \int_{\mathcal{A}_2} \int_{-\infty}^t \begin{bmatrix}
			1 & 0
		\end{bmatrix} e^{A(t-s)} \mathbf{e} \, \Lambda(dA, ds).
\end{equation*}
The eigenvalues of the matrix  
\begin{equation*}
A=\begin{bmatrix}
0 & 1 \\ -a_2 & -a_1
\end{bmatrix}\in  \mathcal{A}_2
\end{equation*}
are given by 
\begin{equation}\label{e:eigenvalues}
    \xi_1=\frac{-a_1 - \sqrt{a_1^2-4a_2}}{2}, \quad \xi_2=\frac{-a_1 + \sqrt{a_1^2-4a_2}}{2}.
\end{equation}
Hence, $A \in \mathcal{A}_2$ if and only if $a_1,a_2>0$. The randomness in $A \in \mathcal{A}_2$ in the definition comes from only two elements, $a_1$ and $a_2$. Thus $\pi$ can be identified with the distribution on $(0,\infty)^2$. The following proposition follows by comparing characteristic functions of the finite dimensional distributions.

\begin{proposition}
Let $\Lambda$ be an $\R$-valued L\'evy basis on $\A_2 \times \R$, $X$ a supCARMA process as in Definition \ref{def:supCARMA} and define $T \colon (0,\infty)^2 \to \A_2$ by
\begin{equation*}
T(a_1,a_2) = \begin{bmatrix}
0 & 1 \\ -a_2 & -a_1
\end{bmatrix}.
\end{equation*}
\begin{enumerate}[(i)]
\item Then $\Lambda'(da_1, da_2, ds)=\Lambda (T(da_1,da_2), ds)$ is an  $\R$-valued L\'evy basis on $(0,\infty)^2 \times \R$ determined by the control measure $\pi'\times \Leb$, $\pi'=\pi \circ T$,  and L\'evy-Khinchine triplet $(a, \Sigma, \nu)$.
\item The process
\begin{equation*}
X'(t) = \int_{(0,\infty)^2} \int_{-\infty}^t \begin{bmatrix}
			1 & 0
		\end{bmatrix} e^{(t-s) 
		\begin{bmatrix}
			0 & 1 \\ -a_2 & -a_1
		\end{bmatrix}
		}  \bm{e} \Lambda'(da_1, da_2, ds), \quad t\in \R,
\end{equation*}
is a version of $X$.
\end{enumerate}
\end{proposition}

We therefore adopt the following definition.

\begin{definition}
Let $\Lambda$ be an $\R$-valued L\'evy basis on $(0,\infty)^2 \times \R$ with generating quadruple $(\gamma,\Sigma,\mu,\pi)$, where $\pi$ is a probability measure on $(0,\infty)^2$. A supCAR$(2)$ process $X$ is defined by
\begin{equation}\label{e:def:supCAR2}
X(t) = \int_{(0,\infty)^2} \int_{\R} g(a_1,a_2,t-s) \Lambda(da_1, da_2, ds), \quad t \in \R,
\end{equation}
where
\begin{equation}\label{e:def:supCAR2:g}
g(a_1,a_2,u) = \begin{bmatrix}
			1 & 0
		\end{bmatrix} e^{u 
		\begin{bmatrix}
			0 & 1 \\ -a_2 & -a_1
		\end{bmatrix}
		}   \begin{bmatrix}
0\\ 1
\end{bmatrix}\1_{(0,\infty)}(u).
\end{equation}
\end{definition}

The matrix exponential function in \eqref{e:def:supCAR2:g} can be computed using the Jordan decomposition by expressing the matrix $A$ as the product $PJP^{-1}$, where $P$ is a regular matrix and $J$ is a Jordan matrix. To calculate the matrices $P$ and $J$, eigenvalues and eigenvectors of the matrix $A$ are required. We have three different cases:
\begin{enumerate}[(I)]
    \item $a_1^2-4a_2=0 \iff \xi_1=\xi_2=-\frac{a_1}{2}$ and we have that 
	\begin{equation*}
	P=\begin{bmatrix}
    1 & 0 \\ 
    -\frac{a_1}{2} & 1
    \end{bmatrix}, \quad
    J=\begin{bmatrix}
    -\frac{a_1}{2} & 1 \\ 
    0 & -\frac{a_1}{2}
    \end{bmatrix}.
	\end{equation*}
	Hence
	\begin{equation*}
	\begin{bmatrix}
		1 & 0
	\end{bmatrix}
	e^{(t-s)A}\mathbf{e}=\begin{bmatrix}
		1 & 0
	\end{bmatrix}Pe^{(t-s)J}P^{-1}\mathbf{e}=(t-s)e^{-\frac{a_1}{2}(t-s)}.    
	\end{equation*}
	
	\item $a_1^2-4a_2>0 \iff \xi_{1,2} = \frac{-a_1 \pm \sqrt{a_1^2-4a_2}}{2}$ and we have that
	\begin{equation*}
		P=\begin{bmatrix}
			1 & 1 \\ 
			\xi_1 & \xi_2
		\end{bmatrix}, \quad
		J=\begin{bmatrix}
			\xi_1 & 0 \\ 
			0 & \xi_2
		\end{bmatrix}.
	\end{equation*}
	Hence
	\begin{equation*}
	\begin{bmatrix}
		1 & 0
	\end{bmatrix}
	e^{(t-s)A}\mathbf{e}=\begin{bmatrix}
		1 & 0
	\end{bmatrix}Pe^{(t-s)J}P^{-1}\mathbf{e}=\frac{e^{\xi_1(t-s)}-e^{\xi_2(t-s)}}{\xi_1 - \xi_2}.
	\end{equation*}

    \item $a_1^2-4a_2<0 \iff \xi_{1,2} = \frac{-a_1 \pm i\sqrt{4a_2-a_1^2}}{2}$ and the matrices $P$ and $J$ have the same form as in the previous case but here $\xi_2=\overline{\xi_1}$.
\end{enumerate}

Let 
\begin{align*}
	\D_{\I}=\{(a_1,a_2)\in (0,\infty)^2\colon a_1^2-4a_2=0\}, \\
	\D_{\II}=\{(a_1,a_2)\in (0,\infty)^2\colon a_1^2-4a_2>0\}, \\
	\D_{\III}=\{(a_1,a_2)\in (0,\infty)^2\colon a_1^2-4a_2<0\},
\end{align*}
and
\begin{equation}\label{e:gfunctions}
\begin{aligned}
	g_{\I}(a_1,a_2,u) &= ue^{-\frac{a_1}{2} u}\1_{[0,\infty)}(u)\1_{\D_{\I}}(a_1,a_2), \\
	g_{\II}(a_1,a_2,u) &= \frac{e^{\xi_1 u}-e^{\xi_2 u}}{\xi_1 - \xi_2} \1_{[0,\infty)}(u)\1_{\mathcal{D}_{\II}}(a_1,a_2)\\	
	&=\frac{1}{\sqrt{\frac{a_1^2}{4}-a_2}}\sinh \left(u\sqrt{\frac{a_1^2}{4}-a_2}\right)e^{-\frac{a_1}{2} u}\1_{[0,\infty)}(u)\1_{\D_{\II}}(a_1,a_2), \\
	g_{\III}(a_1,a_2,u) &= \frac{e^{\xi_1 u}-e^{\xi_2 u}}{\xi_1 - \xi_2} \1_{[0,\infty)}(u)\1_{\mathcal{D}_{\III}}(a_1,a_2)\\	
	&=\frac{1}{\sqrt{a_2-\frac{a_1^2}{4}}} \sin \left(u\sqrt{a_2-\frac{a_1^2}{4}}\right)e^{-\frac{a_1}{2} u}\1_{[0,\infty)}(u)\1_{\D_{\III}}(a_1,a_2), \\
\end{aligned}
\end{equation}
where $\xi_{1,2} = \frac{-a_1 \pm \sqrt{a_1^2-4a_2}}{2}$. Then \eqref{e:def:supCAR2} can be rewritten as a sum of independent mixed moving average processes
\begin{equation*}
	\begin{aligned}
	X(t) &=  X_{\I}(t) + X_{\II} (t) + X_{\III} (t)\\
	&:= \int_{(0,\infty)^2\times \R} g_{\I}(a_1,a_2,t-s)  \, \Lambda(da_1, da_2, ds) + \int_{(0,\infty)^2\times \R} g_{\II}(a_1,a_2,t-s)  \, \Lambda(da_1, da_2, ds) \\
    &\quad+ \int_{(0,\infty)^2\times \R} g_{\III}(a_1,a_2,t-s)  \, \Lambda(da_1, da_2, ds).
	\end{aligned}
\end{equation*}
This implies that we will have essentially three different types of supCAR$(2)$ processes depending on whether $\pi$ is supported in $\D_{\I}$, $\D_{\II}$ or $\D_{\III}$. If $\pi$ is supported on more than one of these sets, then the process can be decomposed into independent components. Hence, it is enough to consider each of the three types separately. We do this in next sections and analyze the properties of these processes.

\section{SupCAR$(2)$ of type I}\label{sec:typeI}

The supCAR$(2)$ of type I corresponds to kernel $g_I$ in \eqref{e:gfunctions}, when the eigenvalues of matrix $A$ are real and equal. In this case, $a_1$ and $a_2$ are directly related. We therefore make the following definition.

\begin{definition}
Let $\rho_{\I}$ be a probability measure on $(0,\infty)$ and $\Lambda_{\I}$ a L\'evy basis on $(0,\infty) \times \R$ corresponding to the L\'evy-Khintchine triplet $(\gamma, \Sigma, \mu)$ and control measure $\rho_{\I}\times \Leb$. Assume that $\int_{|x|>1} \log |x| \, \mu(dx) < \infty$ and
\begin{equation}\label{eq:typeI:cond}
    \int_{(0,\infty)} a^{-3} \, \rho_{\I}(da)<\infty.
\end{equation}
A \textbf{supCAR$(2)$ process of type I} (supCAR$(2)$-I) is defined by
\begin{equation*}
X_{\I}(t) = \int_{(0,\infty)\times \R} (t-s)e^{-\frac{a}{2} (t-s)}\1_{[0,\infty)}(t-s) \, \Lambda_{\I}(da, ds).
\end{equation*}
\end{definition}

If $\rho$ is degenerate, then supCAR$(2)$-I process is a special case of fractional L\'evy-driven OU (fOU) process defined in \cite{wolpert2005}, obtained by taking $\kappa=1$ in the definition given in \cite{wolpert2005}. Hence, one can consider supCAR$(2)$-I as a superposition of fOU processes. In \cite{bngamma2016}, the kernel of supCAR$(2)$-I is called the gamma kernel.

The properties of supCAR$(2)$-I process are determined by L\'evy-Khintchine triplet $(\gamma, \Sigma, \mu)$ and probability measure $\rho_{\I}$. In the following, we will denote by $L$ the L\'evy process with L\'evy-Khintchine triplet $(\gamma, \Sigma, \mu)$.

\begin{theorem}
A supCAR$(2)$-I process is well-defined and stationary with cumulant function of its finite-dimensional distributions given by
    \begin{equation}\label{eq:cum:typeI}
        \kappa_{\left(X_{\I}(t_1), \dots, X_{\I}(t_m)\right)}(\zeta_1, \dots \zeta_m) = \int_{(0,\infty)} \int_{\R} \kappa_{L} \left( \sum_{j=1}^m \zeta_j (t_j-s)e^{-\frac{a}{2} (t_j-s)}\1_{[0,\infty)}(t_j-s) \right) ds \, \rho_{\I}(da),
    \end{equation}
where $t_1 < \cdots < t_m$. The distribution of $X_{\I}(t)$ is infinitely divisible with characteristic triplet $(\gamma_{\I}, \Sigma_{\I}, \mu_{\I})$ given by
\begin{align}
    \gamma_{\I} &= 4 \gamma \int_{(0,\infty)} a^{-2} \, \rho_{\I}(da) + \int_{(0,\infty)} \int_0^{\infty} \int_{\R} \left( \tau(xse^{-\frac{as}{2}})-se^{-\frac{as}{2}}\tau(x) \right) \, \mu(dx) ds \, \rho_{\I}(da),\nonumber\\
    \Sigma_{\I} &= 2 \Sigma \int_{(0,\infty)} a^{-3} \, \rho_{\I}(da),\nonumber\\
    \mu_{\I}(B) &= \int_{(0,\infty)} \int_0^{\infty} \int_{\R} \1_{B}\left(s e^{-\frac{as}{2}} x\right) \, \mu(dx) \, ds \, \rho_{\I}(da), \quad B \in \mathcal{B}(\R).\label{e:muI}
\end{align}
\end{theorem}

\begin{proof}
To show that supCAR$(2)$-I process is well-defined, we will check conditions \eqref{eq:rajros:mma}, which in this case read as
\begin{align}
& \int_{(0,\infty)} \int_0^{\infty} \left| \gamma se^{-\frac{as}{2}} + \int_{\R} \left( \tau(xse^{-\frac{as}{2}})-se^{-\frac{as}{2}}\tau(x) \right) \, \mu(dx) \right|ds \, \rho_{\I}(da) < \infty, \label{eq:condI:1}\\
& \int_{(0,\infty)} \int_0^{\infty} s^2 e^{-as}  \, ds \, \rho_{\I}(da) < \infty, \label{eq:condI:2}\\
& \int_{(0,\infty)} \int_0^{\infty} \int_{\R} \left( 1 \wedge (x^2 s^2 e^{-as})\right) \, \mu(dx) \, ds \, \rho_{\I}(da) < \infty. \label{eq:condI:3}
\end{align}

For \eqref{eq:condI:1}, first note that by \eqref{eq:typeI:cond}
\begin{equation*}
    \int_{(0,\infty)} \int_0^{\infty} se^{-\frac{as}{2}} \, ds \, \rho_{\I}(da) = 4 \int_{(0,\infty)} a^{-2} \, \rho_{\I}(da) <\infty.
\end{equation*}
For the second part we have
\begin{align*}
&\int_{(0,\infty)} \int_0^{\infty} \int_{\R} \left| \tau(xse^{-\frac{as}{2}})-se^{-\frac{as}{2}}\tau(x) \right| \, \mu(dx) \, ds \, \rho_{\I}(da)\\
&\quad = \int_{(0,\infty)} \int_0^{\infty} \int_{|x|\leq 1} |x| se^{-\frac{as}{2}} \1(|x| se^{-\frac{as}{2}} > 1) \, \mu(dx) \, ds \, \rho_{\I}(da)\\
&\quad \quad + \int_{(0,\infty)} \int_0^{\infty} \int_{|x|> 1} |x| se^{-\frac{as}{2}} \1(|x| se^{-\frac{as}{2}} \leq 1)\, \mu(dx) \, ds \, \rho_{\I}(da)\\
&\quad =: J_1+J_2.
\end{align*}
By \eqref{eq:typeI:cond} and since $\mu$ is a L\'evy measure, we get
\begin{align*}
J_1 &\leq \int_{(0,\infty)} \int_0^{\infty} \int_{|x|\leq 1} x^2 s^2 e^{-as} \1(|x| se^{-\frac{as}{2}} > 1) \, \mu(dx) \, ds \, \rho_{\I}(da)\\
&\leq 2 \int_{(0,\infty)} a^{-3} \, \rho_{\I}(da) \int_{|x|\leq 1} x^2 \, \mu(dx) < \infty.
\end{align*}
To show that $J_2<\infty$ we split the integral as follows
\begin{align*}
J_2 &= \int_{(0,\infty)} \int_0^{\infty} \int_{|x|> 1} |x| se^{-\frac{as}{2}} \1(|x| se^{-\frac{as}{2}} \leq 1) \1(|x|\leq ae/2) \, \mu(dx) \, ds \, \rho_{\I}(da)\\
&\quad + \int_{(0,\infty)} \int_{0}^{\infty} \int_{|x|> 1} |x| se^{-\frac{as}{2}} \1(|x| se^{-\frac{as}{2}} \leq 1) \1(|x|> ae/2) \, \mu(dx) \, ds \, \rho_{\I}(da)\\
&=: J_{2,1} + J_{2,2}.
\end{align*}
Consider for the moment the inequality
\begin{equation}\label{eq:acseI:proof}
    |x| se^{-\frac{as}{2}} \leq 1 \iff -\frac{as}{2} e^{-\frac{as}{2}} \geq -\frac{a}{2|x|}.
\end{equation}
Note that the function $h(z)=z e^z$ is decreasing on $(-\infty,-1]$, increasing on $(-1,\infty)$ and $h(-1)=-e^{-1}$. Hence, if $|x|\leq ae/2$, then \eqref{eq:acseI:proof} is valid for any $s>0$. Therefore,
\begin{align*}
J_{2,1} &= \int_{(0,\infty)} \int_0^{\infty} \int_{|x|> 1} |x| se^{-\frac{as}{2}} \1(|x|\leq ae/2) \, \mu(dx) \, ds \, \rho_{\I}(da)\\
&= 4 \int_{(0,\infty)} a^{-2} \int_{|x|> 1} |x| \1(|x|\leq ae/2) \, \mu(dx) \, \rho_{\I}(da)\\
&\leq 2e \int_{(0,\infty)} a^{-1} \1(ae/2>1) \, \rho_{\I}(da) \int_{|x|> 1}  \, \mu(dx) < \infty,
\end{align*}
while for $J_{2,2}$
\begin{align*}
J_{2,2} &= \int_{(0,\infty)} \int_{0}^{2/a} \int_{|x|> 1} |x| se^{-\frac{as}{2}} \1(|x| se^{-\frac{as}{2}} \leq 1) \1(|x|> ae/2) \, \mu(dx) \, ds \, \rho_{\I}(da)\\
&\quad + \int_{(0,\infty)} \int_{2/a}^\infty \int_{|x|> 1} |x| se^{-\frac{as}{2}} \1(|x| se^{-\frac{as}{2}} \leq 1) \1(|x|> ae/2) \, \mu(dx) \, ds \, \rho_{\I}(da)\\
&\leq \int_{(0,\infty)} \int_{0}^{2/a} \int_{|x|> 1} |x| se^{-\frac{as}{2}} \1(|x| se^{-1} \leq 1) \1(|x|> ae/2) \, \mu(dx) \, ds \, \rho_{\I}(da)\\
&\quad + \int_{(0,\infty)} \int_{2/a}^\infty \int_{|x|> 1} |x| se^{-\frac{as}{2}} \1(|x| (2/a) e^{-\frac{as}{2}} \leq 1) \1(|x|> ae/2) \, \mu(dx) \, ds \, \rho_{\I}(da)\\
&\leq \int_{(0,\infty)} \int_{0}^{2/a} \int_{|x|> 1} e \, \mu(dx) \, ds \, \rho_{\I}(da)\\
&\quad + \int_{(0,\infty)} \int_{2/a}^\infty \int_{|x|> 1} |x| se^{-\frac{as}{2}} \1(s \geq 2 \log (2|x|/a)/a) \1(|x|> ae/2) \, \mu(dx) \, ds \, \rho_{\I}(da)\\
&\leq 2 e \int_{(0,\infty)} a^{-1} \, \rho_{\I}(da) \int_{|x|> 1} \, \mu(dx) \\
&\quad + 2 \int_{(0,\infty)} \int_{|x|> 1} |x| (a |x|)^{-1} (1+\log (2|x|/a)) \1(|x|> ae/2) \, \mu(dx) \, \rho_{\I}(da)\\
&\leq 2 e \int_{(0,\infty)} a^{-1} \, \rho_{\I}(da) \int_{|x|> 1} \, \mu(dx) \\
&\quad + 2(1+\log 2) \int_{(0,\infty)} a^{-1} \, \rho_{\I}(da) \int_{|x|> 1} \, \mu(dx) + 2 \int_{(0,\infty)} |\log a| a^{-1} \, \rho_{\I}(da) \int_{|x|> 1} \, \mu(dx)\\
&\quad \quad + 2 \int_{(0,\infty)} a^{-1} \, \rho_{\I}(da) \int_{|x|> 1} \log |x| \, \mu(dx) < \infty.
\end{align*}

If there is a Gaussian component $\Sigma>0$, then \eqref{eq:condI:2} holds since
\begin{equation*}
\int_{(0,\infty)} \int_0^{\infty} s^2 e^{-as} \, ds \, \rho_{\I}(da) = 2 \int_{(0,\infty)} a^{-3} \, \rho_{\I}(da)<\infty.
\end{equation*}

It remains to check \eqref{eq:condI:3}. To this end, let
\begin{align*}
& \int_{(0,\infty)} \int_0^{\infty} \int_{\R} \left( 1 \wedge (x^2 s^2 e^{-as})\right) \, \mu(dx) \, ds \, \rho_{\I}(da)\\
&\quad =\int_{(0,\infty)} \int_0^{\infty} \int_{|x|\leq 1} \left( 1 \wedge (x^2 s^2 e^{-as})\right) \, \mu(dx) \, ds \, \rho_{\I}(da) \\
&\quad \quad + \int_{(0,\infty)} \int_0^{\infty} \int_{|x|>1} \left( 1 \wedge (x^2 s^2 e^{-as})\right) \, \mu(dx) \, ds \, \rho_{\I}(da) =: I_1 + I_2.
\end{align*}
For $I_1$, it follows that
\begin{align*}
I_1 &\leq \int_{(0,\infty)} \int_0^{\infty} \int_{|x|\leq 1} x^2 s^2 e^{-as} \, \mu(dx) \, ds \, \rho_{\I}(da)\\
&= 2 \int_{(0,\infty)} a^{-3} \, \rho_{\I}(da) \int_{|x|\leq 1} x^2 \, \mu(dx) < \infty.
\end{align*}
Next, we split $I_2$ as
\begin{align*}
I_2 & = \int_{(0,\infty)} \int_0^{\infty} \int_{|x|>1} x^2 s^2 e^{-as} \1(|x| se^{-\frac{as}{2}} \leq 1) \, \mu(dx) \, ds \, \rho_{\I}(da)\\
&\quad + \int_{(0,\infty)} \int_0^{\infty} \int_{|x|>1} \1(|x| se^{-\frac{as}{2}} > 1) \, \mu(dx) \, ds \, \rho_{\I}(da) =: I_{2,1} + I_{2,2}.
\end{align*}
Note that $I_{2,1}\leq J_2$, therefore $I_{2,1}<\infty$. For $I_{2,2}$, note that if \eqref{eq:acseI:proof} holds, then we must have $|x|>ae/2$. If $|x|>ae/2$, then the equation $h(z)=z e^z= -\frac{a}{2|x|}$ has two solutions:
\begin{equation*}
z_1 = W_0(-a/(2|x|)) \in [-1,0), \quad z_2 = W_{-1}(-a/(2|x|))\in (-\infty,-1],
\end{equation*}
where $W_0$ and $W_{-1}$ are Lambert $W$ functions (see e.g.~\cite{mezo2022}). Then it follows that $h(z)\geq -\frac{a}{2|x|}$ for $z \in (-\infty, z_2] \cup [z_1,\infty)$. We conclude that, if $|x|>ae/2$, then \eqref{eq:acseI:proof} holds for $s \in (0,s_1] \cup [s_2, \infty)$, where
\begin{equation*}
s_1 = -\frac{2}{a} W_0(-a/(2|x|)), \quad s_2 = -\frac{2}{a} W_{-1}(-a/(2|x|)).
\end{equation*}
Returning back to $I_{2,2}$, we have the following
\begin{align*}
I_{2,2} &= \int_{(0,\infty)} \int_{s_1}^{s_2} \int_{|x|>1} \1(|x| se^{-\frac{as}{2}} > 1) \1(|x|> ae/2) \, \mu(dx) \, ds \, \rho_{\I}(da)\\
&= 2 \int_{(0,\infty)} \int_{|x|>1} a ^{-1} \left( W_0(-a/(2|x|)) - W_{-1}(-a/(2|x|)) \right) \1(|x|> ae/2) \, \mu(dx) \, ds \, \rho_{\I}(da).
\end{align*}
From \cite[Theorem 1]{chatzigeorgiou2013bounds} we get that
\begin{equation*}
    W_{-1}(-a/(2|x|)) > - \log (2 |x|/a) - \sqrt{2(\log (2 |x|/a) - 1)}.
\end{equation*}
For $y>1$, $\sqrt{2(y-1)}<y$, hence
\begin{align*}
I_{2,2} & \leq 4 \int_{(0,\infty)} \int_{|x|>1} a^{-1} \log (2|x|/a) \1(|x|> ae/2)\, \mu(dx) \, \rho_{\I}(da)\\
& \leq 4 \log 2 \int_{(0,\infty)} a^{-1} \, \rho_{\I}(da) \int_{|x|>1} \, \mu(dx) + 4\int_{(0,\infty)} |\log a| a^{-1} \, \rho_{\I}(da) \int_{|x|>1} \, \mu(dx) \\
&\quad + 4 \int_{(0,\infty)} a^{-1} \, \rho_{\I}(da) \int_{|x|>1}\log |x| \, \mu(dx) <\infty.
\end{align*}
From \cite[Proposition 2.6]{rajput1989spectral} we obtain \eqref{eq:cum:typeI}, from which stationarity follows, while \cite[Theorem 2.7]{rajput1989spectral} yields the infinite divisibility of $X(t)$ with characteristic triplet $(\gamma_{\I}, \Sigma_{\I}, \mu_{\I})$.
\end{proof}

A closer inspection of the proof reveals that it is not necessary to assume that $\rho_I$ is a finite measure. In fact, it suffices that $\int_{(1,\infty)} |\log a| a^{-1} \, \rho_{\I}(da)<\infty$.

The following proposition establishes the conditions for the finiteness of moments of $X_{\I}$.

\begin{proposition}\label{prop:mc:I}
    \begin{enumerate}[(i)]
        \item If $0 < p \leq 2$ and
            \begin{equation*}
                \int_{|x|>1} |x|^p \, \mu(dx) < \infty,
            \end{equation*}
            then $\E |X_{\I}(0)|^p < \infty$.
        \item If $p > 2$ and
            \begin{equation*}
                \int_{|x|>1} |x|^p \, \mu(dx) < \infty, \quad \int_{(0,\infty)} a^{-(p+1)} \, \rho_{\I}(da) < \infty,
            \end{equation*}
            then $\E |X_{\I}(0)|^p < \infty$.
    \end{enumerate}
\end{proposition}

\begin{proof}
    According to \cite[Corollary~25.8]{Sato}, it suffices to show that $\int_{|x|>1} |x|^p \, \mu_{\I}(dx)<\infty$ for $\mu_I$ given by \eqref{e:muI}. For $0 < p \leq 2$ we obtain
    \begin{align*}
        \int_{|x|>1} |x|^p \, \mu_{\I}(dx) &= \int_{(0,\infty)} \int_0^{\infty} \int_{\R} s^p e^{-\frac{pas}{2}} |x|^p \1 (|x| s e^{-\frac{as}{2}}>1 ) \, \mu(dx) \, ds \, \rho_{\I}(da) \\
        &= \int_{(0,\infty)} \int_0^{\infty} \int_{|x| \leq 1} s^p e^{-\frac{pas}{2}} |x|^p \1 (|x| s e^{-\frac{as}{2}}>1 ) \, \mu(dx) \, ds \, \rho_{\I}(da) \\
        & \quad + \int_{(0,\infty)} \int_0^{\infty} \int_{|x| > 1} s^p e^{-\frac{pas}{2}} |x|^p \1 (|x| s e^{-\frac{as}{2}}>1 ) \, \mu(dx) \, ds \, \rho_{\I}(da) \\
        & \leq \int_{(0,\infty)} \int_0^{\infty} \int_{|x| \leq 1} s^2 e^{-as} |x|^2 \1 (|x| s e^{-\frac{as}{2}}>1 ) \, \mu(dx) \, ds \, \rho_{\I}(da) \\
        & \quad + \int_{(0,\infty)} \int_0^{\infty} \int_{|x| > 1} s^p e^{-\frac{pas}{2}} |x|^p \1 (|x| s e^{-\frac{as}{2}}>1 ) \, \mu(dx) \, ds \, \rho_{\I}(da) \\
        & \leq 2 \int_{(0,\infty)} a^{-3} \, \rho_{\I}(da) \int_{|x| \leq 1} |x|^2 \, \mu(dx) \\
        & \quad + \left( \frac{2}{p} \right)^{p+1} \Gamma(p+1) \int_{(0,\infty)} a^{-(p+1)} \, \rho_{\I}(da) \int_{|x| > 1} |x|^p \, \mu(dx) < \infty.
    \end{align*}
    Since \eqref{eq:typeI:cond} implies $\int_{(0,\infty)} a^{-(p+1)} \rho_{\I}(da) < \infty$, and the remaining terms are finite by assumption and the fact that $\mu$ is a L\'evy measure, the above expression is finite.

    For $p>2$, we already have
    \begin{equation*}
        \int_{\R} |x|^p \, \mu(dx) = \int_{|x|\leq 1} |x|^p \, \mu(dx) + \int_{|x|>1} |x|^p \, \mu(dx) \leq \int_{|x|\leq 1} |x|^2 \, \mu(dx) + \int_{|x|>1} |x|^p \, \mu(dx) < \infty.
    \end{equation*}
    Therefore,
    \begin{align*}
        \int_{|x|>1} |x|^p \, \mu_{\I}(dx) &= \int_{(0,\infty)} \int_0^{\infty} \int_{\R} s^p e^{-\frac{pas}{2}} |x|^p \1 (|x| s e^{-\frac{as}{2}}>1 ) \, \mu(dx) \, ds \, \rho_{\I}(da) \\
        & \leq \int_{(0,\infty)} \int_0^{\infty} s^p e^{-\frac{pas}{2}} \, ds \, \rho_{\I}(da) \int_{\R} |x|^p \, \mu(dx) \\
        & \leq \left( \frac{2}{p} \right)^{p+1} \Gamma(p+1) \int_{(0,\infty)} a^{-(p+1)} \, \rho_{\I}(da) \int_{\R} |x|^p \, \mu(dx) < \infty.
    \end{align*}
\end{proof}

\begin{proposition}\label{prop:corr:typeI}
    If $\int_{|x|>1} x^2 \, \mu(dx) < \infty$, then the process $X_{\I}$ has finite second moment and
    \begin{align*}
        \E(X_{\I}(t)) &= 4 m_{-2}(\rho_{\I}) \E(L(1)),\\
        \Var(X_{\I}(t)) &= 2 m_{-3}(\rho_{\I}) \Var(L(1)),
    \end{align*}
    while the correlation function is
    \begin{equation*}
        r(\tau) = \frac{1}{m_{-3}(\rho_{\I})} \int_{(0,\infty)} a^{-3} \left(\frac{a\tau}{2}+1\right) e^{-\frac{a\tau}{2}} \rho_{\I}(da).
    \end{equation*}
\end{proposition}

\begin{proof}
    The finiteness of the second moment follows from Proposition \ref{prop:mc:I}. Differentiating \eqref{eq:cum:typeI} with respect to the appropriate arguments and evaluating at zero yields the desired cumulants. In fact, for $m=1$ in \eqref{eq:cum:typeI}, differentiating with respect to $\zeta_1$ gives
    \begin{equation*}
        \kappa_{X_{\I}}'(0) = \kappa_{L}'(0) \int_{(0,\infty)} \int_{-\infty}^{t_1} (t_1-s)e^{-\frac{a}{2} (t_1-s)} \, ds \, \rho_{\I}(da).
    \end{equation*}
    Since $\E(X_{\I}(t))=\kappa_{X_{\I}}^{(1)}=-i\kappa_{X_{\I}}'(0)$, we get
    \begin{equation*}
        \E(X_{\I}(t)) = 4 \E(L(1)) \int_{(0,\infty)} a^{-2} \, \rho_{\I}(da) = 4 m_{-2}(\rho_{\I}) \E(L(1)).
    \end{equation*}
    For $m=2$, differentiating with respect to $\zeta_1$ and $\zeta_2$, and evaluating at $(\zeta_1, \zeta_2) = (0,0)$ gives
    \begin{equation*}
        \Cov(X_{\I}(t_1), X_{\I}(t_2)) = - \kappa_L''(0) \int_{(0,\infty)} \int_{\R} (t_1-s)(t_2-s) e^{-\frac{a}{2} (t_1+t_2-2s)}\1_{[0,\infty)}(t_1-s)\1_{[0,\infty)}(t_2-s) \, ds \, \rho_{\I}(da).
    \end{equation*}
    Using $-\kappa_L''(0) = \Var(L(1))$, integration with respect to $s$ gives
    \begin{equation*}
        \Cov(X_{\I}(t_1), X_{\I}(t_2)) = 2 \Var(L(1)) \int_{(0,\infty)} a^{-3} \left(\frac{a(t_2-t_1)}{2}+1\right) e^{-\frac{a}{2}(t_2-t_1)} \rho_{\I}(da).
    \end{equation*}
    By stationarity 
        $$\Cov(X_{\I}(t_1), X_{\I}(t_2)) = \Cov(X_{\I}(0), X_{\I}(t_2-t_1))$$
    and for $\tau>0$, we obtain
    \begin{equation*}
        \Cov(X_{\I}(0), X_{\I}(\tau)) = 2 \Var(L(1)) \int_{(0,\infty)} a^{-3} \left(\frac{a\tau}{2}+1\right) e^{-\frac{a\tau}{2}} \, \rho_{\I}(da).
    \end{equation*}
    In particular, the variance is
    \begin{equation*}
        \Var(X_{\I}(0)) = 2 \Var(L(1)) \int_{(0,\infty)} a^{-3} \, \rho_{\I}(da) = 2 {m_{-3}(\rho_{\I})} \Var(L(1)),
    \end{equation*}
    hence the correlation function is
    \begin{equation*}
        r(\tau)=\frac{1}{{m_{-3}(\rho_{\I})}} \int_{(0,\infty)} a^{-3} \left(\frac{a\tau}{2}+1\right) e^{-\frac{a\tau}{2}} \rho_{\I}(da).
    \end{equation*}
\end{proof}

\begin{example}
	Suppose that $\rho_{\I}$ is a $\Gamma(\alpha+3,1)$ distribution given by
    \begin{equation*}
        \rho_{\I}(da) = \frac{1}{\Gamma(\alpha+3)} a^{\alpha+2}e^{-a}\,da,
    \end{equation*}
    where $\alpha>0$. In particular, condition \eqref{eq:typeI:cond} is satisfied. Assuming $\int_{|x|>1} x^2 \, \mu(dx) < \infty$, it follows that
	\begin{align*}
		{m_{-3}(\rho_{\I})}&=\int_{(0,\infty)} a^{-3} \, \rho_{\I}(da)=\frac{1}{\alpha(\alpha+1)(\alpha+2)}
	\end{align*}
	and
	\begin{align*}
		r(\tau) &=\frac{1}{{m_{-3}(\rho_{\I})}}\int_0^{\infty} a^{-3} \left(\frac{a\tau}{2}+1\right) e^{-\frac{a\tau}{2}} \frac{1}{\Gamma(\alpha+3)} a^{\alpha+2}e^{-a}\,da\\
		&=\frac{\alpha}{2} \left( \tau\left( \frac{\tau}{2}+1\right)^{-(\alpha+1)}+\frac{2}{\alpha} \left( \frac{\tau}{2}+1\right)^{-\alpha}\right) \\
		&=\frac{2^{\alpha}}{(\tau+2)^{\alpha+1}}(\alpha \tau+\tau+2).
	\end{align*}
	Therefore, the correlation function satisfies
	\begin{equation*}
	    r(\tau)\sim 2^{\alpha}(\alpha+1)\tau^{-\alpha},
	\end{equation*}
	and for $\alpha\in(0,1]$ the correlation function exhibits the property of long-range dependence.
\end{example}

\section{SupCAR$(2)$ of type II}\label{sec:typeII}

We now consider the case when the eigenvalues of matrix $A$ are real and distinct, hence the kernel function is $g_{\II}$ from \eqref{e:gfunctions}. Since the relation between matrix elements $a_1, a_2$ and the kernel function is complicated, it is more natural and convenient to define the process by specifying the distribution of eigenvalues. To this end, note that since the eigenvalues \eqref{e:eigenvalues} are real and distinct, they can always be written as 
\begin{equation*}
	\xi_1=-\lambda, \quad \xi_2=-\lambda \theta,
\end{equation*}
where $\lambda>0$ and $\theta \in (0,1)$. We introduce randomness in eigenvalues through $\lambda$ and $\theta$. We will show that this allows tailoring the distributional properties and dependence structure by specifying the distribution of $\lambda$ and $\theta$.

\begin{definition}
Let $\pi_{\lambda, \theta}$ be a probability measure on $(0,\infty) \times (0,1)$ and $\Lambda_{\II}$ a L\'evy basis on $(0,\infty) \times (0,1) \times \R$ corresponding to the L\'evy-Khintchine triplet $(\gamma, \Sigma, \mu)$ and control measure $\pi_{\lambda, \theta} \times \Leb$. Assume that $\int_{|x|>1} \log |x| \, \mu(dx)$ and
\begin{equation}\label{eq:typeII:cond}
\begin{aligned}
    &\int_{(0,\infty) \times (0,1)} \frac{1}{\lambda^3 \theta} \, \pi_{\lambda, \theta}(d\lambda, d\theta) < \infty, \\ &\int_{(0,\infty) \times (0,1)} \frac{|\log(1-\theta)|}{\lambda} \, \pi_{\lambda, \theta}(d\lambda, d\theta) < \infty, \\
    &\int_{(0,1)} \frac{1}{\theta} \, \pi_{\theta}(d\theta) < \infty.
\end{aligned}
\end{equation}
A \textbf{supCAR$(2)$ process of type II} (supCAR$(2)$-II) is defined by
\begin{equation*}
    X_{\II}(t) = \int_{(0,\infty) \times (0,1) \times \R} \frac{e^{-\lambda \theta (t-s)}-e^{-\lambda  (t-s)}}{\lambda(1-\theta)} \1_{[0,\infty)}(t-s)\, \Lambda_{\II}(d\lambda, d\theta, ds).
\end{equation*}
\end{definition}

\begin{theorem}
A supCAR$(2)$-II process is well-defined and stationary with cumulant function of its finite-dimensional distributions given by
    \begin{align}\label{eq:cum:typeII}
        &\kappa_{\left(X_{\II}(t_1), \dots, X_{\II}(t_m)\right)}(\zeta_1, \dots \zeta_m) \notag \\
        &\quad = \int_{(0,\infty) \times (0,1)} \int_{\R} \kappa_{L} \left( \sum_{j=1}^m \zeta_j \frac{e^{-\lambda \theta (t_j-s)}-e^{-\lambda  (t_j-s)}}{\lambda(1-\theta) } \1_{[0,\infty)}(t_j-s) \right) ds \, \pi_{\lambda, \theta}(d\lambda, d\theta),
    \end{align}
where $t_1 < \cdots < t_m$. The distribution of $X_{\II}(t)$ is infinitely divisible with characteristic triplet $(\gamma_{\II}, \Sigma_{\II}, \mu_{\II})$ given by
\begin{align*}
    \gamma_{\II} &= \gamma \int_{(0,\infty) \times (0,1)} \frac{1}{\lambda^2 \theta} \, \pi_{\lambda, \theta}(d\lambda, d\theta) \\
    & \quad + \int_{(0,\infty) \times (0,1)} \int_0^{\infty} \int_{\R} \left( \tau\left(x\frac{e^{-\lambda \theta s}-e^{-\lambda s}}{\lambda(1-\theta)}\right)-\frac{e^{-\lambda \theta s}-e^{-\lambda s}}{\lambda(1-\theta)}\tau(x) \right) \, \mu(dx) \, ds \, \pi_{\lambda, \theta}(d\lambda, d\theta),\\
    \Sigma_{\II} &= \frac{\Sigma}{2} \int_{(0,\infty) \times (0,1)} \frac{1}{\lambda^3 \theta(1+\theta)} \, \pi_{\lambda, \theta}(d\lambda, d\theta), \\
    \mu_{\II}(B) &= \int_{(0,\infty) \times (0,1)} \int_0^{\infty} \int_{\R} \1_{B}\left(\frac{e^{-\lambda \theta s}-e^{-\lambda s}}{\lambda(1-\theta)} x\right) \, \mu(dx) \, ds \, \pi_{\lambda, \theta}(d\lambda, d\theta), \quad B \in \mathcal{B}(\R).
\end{align*}
\end{theorem}

\begin{proof}
    Let
    \begin{equation*}
        f(\lambda, \theta, s) := \frac{e^{-\lambda \theta s}-e^{-\lambda s}}{\lambda(1-\theta)} \1_{[0,\infty)}(s).
    \end{equation*}
    We check that the integrability requirements \eqref{eq:rajros:mma} are satisfied. In this case, the conditions can be written as
    \begin{align}
        & \int_{(0,\infty) \times (0,1)} \int_{\R} \Big| \gamma f(\lambda, \theta, s) + \int_{\R} \left( \tau(xf(\lambda, \theta, s))-f(\lambda, \theta, s)\tau(x) \right) \, \mu(dx) \Big| \, ds \, \pi_{\lambda, \theta}(d\lambda, d\theta) < \infty, \label{eq:condII:1}\\
        & \int_{(0,\infty) \times (0,1)} \int_{\R} f(\lambda, \theta, s)^2  \, ds \, \pi_{\lambda, \theta}(d\lambda, d\theta) < \infty, \label{eq:condII:2}\\
        & \int_{(0,\infty) \times (0,1)} \int_{\R} \int_{\R} \left( 1 \wedge (x^2 f(\lambda, \theta, s)^2)\right) \, \mu(dx) \, ds \, \pi_{\lambda, \theta}(d\lambda, d\theta) < \infty. \label{eq:condII:3}
    \end{align}
    First, note that since $1/\lambda \leq 1/\lambda^3$ for $\lambda \in (0,1)$ and $1/\lambda \leq 1$ for $\lambda\geq1$, we have
    \begin{equation}\label{e:proofII:1}
    \begin{aligned}
        \int_{(0,\infty) \times (0,1)} \frac{1}{\lambda \theta} \, \pi_{\lambda, \theta}(d\lambda, d\theta) &= \int_{(0,1) \times (0,1)} \frac{1}{\lambda \theta} \, \pi_{\lambda, \theta}(d\lambda, d\theta) + \int_{[1,\infty) \times (0,1)} \frac{1}{\lambda \theta} \, \pi_{\lambda, \theta}(d\lambda, d\theta)\\
        & \leq \int_{(0,1) \times (0,1)} \frac{1}{\lambda^3 \theta} \, \pi_{\lambda, \theta}(d\lambda, d\theta) + \int_{[1,\infty) \times (0,1)} \frac{1}{\theta} \, \pi_{\lambda, \theta}(d\lambda, d\theta) \\
        & \leq \int_{(0,1) \times (0,1)} \frac{1}{\lambda^3 \theta} \, \pi_{\lambda, \theta}(d\lambda, d\theta) + \int_{(0,1)} \frac{1}{\theta} \, \pi_{\theta}(d\theta) < \infty.
    \end{aligned}
    \end{equation}
    Similarly, since $1/\lambda^{2} \leq 1/\lambda^{3}$ for $\lambda \in (0,1)$ and $1/\lambda^{2} \leq 1/\lambda$ for $\lambda \geq 1$, it follows that
    \begin{equation}\label{e:proofII:2}
    \begin{aligned}
        \int_{(0,\infty) \times (0,1)} \frac{1}{\lambda^2 \theta} \, \pi_{\lambda, \theta}(d\lambda, d\theta) &= \int_{(0,1) \times (0,1)} \frac{1}{\lambda^2 \theta} \, \pi_{\lambda, \theta}(d\lambda, d\theta) + \int_{[1,\infty) \times (0,1)} \frac{1}{\lambda^2 \theta} \, \pi_{\lambda, \theta}(d\lambda, d\theta)\\
        & \leq \int_{(0,1) \times (0,1)} \frac{1}{\lambda^3 \theta} \, \pi_{\lambda, \theta}(d\lambda, d\theta) + \int_{[1,\infty) \times (0,1)} \frac{1}{\lambda \theta} \, \pi_{\lambda, \theta}(d\lambda, d\theta) < \infty.
    \end{aligned}
    \end{equation}
    Furthermore, since the function $x \mapsto x \log(1/x)$ reaches its maximum at $e^{-1}$ with value $e^{-1}$, we have $\log(1/\lambda) \leq e^{-1} \lambda^{-1}$ for all $\lambda>0$, and hence
    \begin{equation}\label{e:proofII:3}
    \begin{aligned}
        &\int_{(0,\infty) \times (0,1)} \frac{\log(1/\lambda)}{\lambda \theta} \, \pi_{\lambda, \theta}(d\lambda, d\theta) \\
        &\quad = \int_{(0,1) \times (0,1)} \frac{\log(1/\lambda)}{\lambda \theta} \, \pi_{\lambda, \theta}(d\lambda, d\theta) + \int_{[1,\infty) \times (0,1)} \frac{\log(1/\lambda)}{\lambda \theta} \, \pi_{\lambda, \theta}(d\lambda, d\theta)\\
        &\quad \leq \frac{1}{e} \int_{[1,\infty) \times (0,1)} \frac{1}{\lambda^2 \theta} \, \pi_{\lambda, \theta}(d\lambda, d\theta) < \infty.
    \end{aligned}
    \end{equation}
    
    To establish \eqref{eq:condII:1}, note first that
    \begin{equation*} 
        \int_{(0,\infty) \times (0,1)} \int_{\R} |f(\lambda, \theta, s)| \, ds \, \pi_{\lambda, \theta}(d\lambda, d\theta) = \int_{(0,\infty) \times (0,1)} \frac{1}{\lambda^2 \theta} \, \pi_{\lambda, \theta}(d\lambda, d\theta),
    \end{equation*}
    which is finite by \eqref{e:proofII:2}. For the second part we have
    \begin{equation}\label{e:proofII:J2}
    \begin{aligned}
        &\int_{(0,\infty) \times (0,1)} \int_{\R} \int_{\R} \left| \tau \left(xf(\lambda, \theta, s)\right)-f(\lambda, \theta, s)\tau(x) \right| \, \mu(dx) \, ds \, \pi_{\lambda, \theta}(d\lambda, d\theta)\\
        &\quad = \int_{(0,\infty) \times (0,1)} \int_{\R} \int_{|x|\leq 1} |x|f(\lambda, \theta, s) \1(|x|f(\lambda, \theta, s)> 1) \, \mu(dx) \, ds \, \pi_{\lambda, \theta}(d\lambda, d\theta)\\
        &\quad \quad + \int_{(0,\infty) \times (0,1)} \int_{\R} \int_{|x|> 1} |x|f(\lambda, \theta, s) \1(|x|f(\lambda, \theta, s) \leq 1)\, \mu(dx) \, ds \, \pi_{\lambda, \theta}(d\lambda, d\theta)\\
        &\quad =: J_1+J_2.
    \end{aligned}
    \end{equation}
        For $J_1$, it follows by \eqref{eq:typeII:cond} that
    \begin{align*}
        J_1 & \leq  \int_{(0,\infty) \times (0,1)} \int_{\R} \int_{|x|\leq 1} x^2 f(\lambda, \theta, s)^2 \1(|x|f(\lambda, \theta, s)> 1) \, \mu(dx) \, ds \, \pi_{\lambda, \theta}(d\lambda, d\theta)\\
        & \leq \int_{(0,\infty) \times (0,1)} \int_0^{\infty} \int_{|x|\leq 1} x^2 \left(\frac{e^{-\lambda \theta s}-e^{-\lambda s}}{\lambda(1-\theta)}\right)^2 \, \mu(dx) \, ds \, \pi_{\lambda, \theta}(d\lambda, d\theta) \\
        & \leq \frac{1}{2} \int_{(0,\infty) \times (0,1)} \frac{1}{\lambda^3 \theta(1+\theta)} \, \pi_{\lambda, \theta}(d\lambda, d\theta) \int_{|x|\leq 1} x^2 \, \mu(dx) \\
        & \leq \frac{1}{2} \int_{(0,\infty) \times (0,1)} \frac{1}{\lambda^3 \theta} \, \pi_{\lambda, \theta}(d\lambda, d\theta) \int_{|x|\leq 1} x^2 \, \mu(dx) < \infty.
    \end{align*}
    Since
    \begin{equation*}
        |x| \frac{e^{-\lambda \theta s}-e^{-\lambda s}}{\lambda(1-\theta)} \leq |x| \frac{e^{-\lambda \theta s}}{\lambda(1-\theta)},
    \end{equation*}
    it follows that
    \begin{equation}\label{eq:proofII:ineq}
        |x| \frac{e^{-\lambda \theta s}-e^{-\lambda s}}{\lambda(1-\theta)} \leq 1
    \end{equation}
    whenever $|x| \frac{e^{-\lambda \theta s}}{\lambda(1-\theta)} \leq 1$. Define
    \begin{equation}\label{e:proofII:s0}
        s_0 := \frac{\log \left(|x| / (\lambda(1-\theta))\right)}{\lambda \theta}.
    \end{equation}
    Then, for $s>s_0$, inequality \eqref{eq:proofII:ineq} certainly holds. If $|x| \leq \lambda(1-\theta)$, we have $s_0<0$, and \eqref{eq:proofII:ineq} holds for all $s>0$. Moreover, $s_0>0$ whenever $|x|>\lambda(1-\theta)$. Consequently, we decompose $J_2$ as
    \begin{align*}
        J_2 & = \int_{(0,\infty) \times (0,1)} \int_{\R} \int_{|x|> 1} |x|f(\lambda, \theta, s) \1(|x|f(\lambda, \theta, s) \leq 1) \1(|x| \leq \lambda(1-\theta))\, \mu(dx) \, ds \, \pi_{\lambda, \theta}(d\lambda, d\theta) \\
        & \quad + \int_{(0,\infty) \times (0,1)} \int_{\R} \int_{|x|> 1} |x|f(\lambda, \theta, s) \1(|x|f(\lambda, \theta, s) \leq 1) \1(|x| > \lambda(1-\theta))\, \mu(dx) \, ds \, \pi_{\lambda, \theta}(d\lambda, d\theta) \\
        & =: J_{2,1} + J_{2,2}.
    \end{align*}
    For $J_{2,1}$, we get
    \begin{align*}
        J_{2,1} & = \int_{(0,\infty) \times (0,1)} \int_{\R} \int_{|x|> 1} |x|f(\lambda, \theta, s) \1(|x| \leq \lambda(1-\theta))\, \mu(dx) \, ds \, \pi_{\lambda, \theta}(d\lambda, d\theta) \\
        & = \int_{(0,\infty) \times (0,1)} \frac{1}{\lambda^2 \theta} \int_{|x|> 1} |x| \1(|x| \leq \lambda(1-\theta))\, \mu(dx) \, \pi_{\lambda, \theta}(d\lambda, d\theta)\\
        & \leq \int_{(0,\infty) \times (0,1)} \frac{1-\theta}{\lambda \theta} \1(\lambda(1-\theta)>1) \, \pi_{\lambda, \theta}(d\lambda, d\theta) \int_{|x|> 1} \, \mu(dx)\\
        & \leq \int_{(0,\infty) \times (0,1)} \frac{1}{\lambda \theta} \, \pi_{\lambda, \theta}(d\lambda, d\theta) \int_{|x|> 1} \, \mu(dx),
    \end{align*}
    and the finiteness follows from \eqref{e:proofII:1}. For $J_{2,2}$ it follows that
    \begin{equation}\label{e:proofII:J22}
    \begin{aligned}
        J_{2,2} & \leq \int_{(0,\infty) \times (0,1)} \int_0^{s_0} \int_{|x|> 1} \1(|x| > \lambda(1-\theta)) \, \mu(dx) \, ds \, \pi_{\lambda, \theta}(d\lambda, d\theta) \\
        & \quad + \int_{(0,\infty) \times (0,1)} \int_{s_0}^{\infty} \int_{|x|> 1} |x|f(\lambda, \theta, s) \, \mu(dx) \, ds \, \pi_{\lambda, \theta}(d\lambda, d\theta)\\
        & \leq \int_{(0,\infty) \times (0,1)} \int_{|x|> 1} \frac{\log(|x| / (\lambda(1-\theta)))}{\lambda \theta} \1(|x| > \lambda(1-\theta)) \, \mu(dx) \, \pi_{\lambda, \theta}(d\lambda, d\theta) \\
        & \quad + \int_{(0,\infty) \times (0,1)} \int_{|x|> 1} \frac{1}{\lambda \theta} \, \mu(dx) \, \pi_{\lambda, \theta}(d\lambda, d\theta) \\
        & \leq \int_{(0,\infty) \times (0,1)} \frac{1}{\lambda \theta} \, \pi_{\lambda, \theta}(d\lambda, d\theta) \int_{|x|> 1} \log |x| \, \mu(dx) \\
        & \quad + \int_{(0,\infty) \times (0,1)} \frac{\log(1/\lambda)}{\lambda \theta} \, \pi_{\lambda, \theta}(d\lambda, d\theta) \int_{|x|> 1} \, \mu(dx) \\
        & \quad \quad + \int_{(0,\infty) \times (0,1)} \frac{|\log(1-\theta)|}{\lambda \theta} \, \pi_{\lambda, \theta}(d\lambda, d\theta) \int_{|x|> 1} \, \mu(dx) \\
        & \quad \quad \quad + \int_{(0,\infty) \times (0,1)} \frac{1}{\lambda \theta} \, \pi_{\lambda, \theta}(d\lambda, d\theta) \int_{|x|> 1} \, \mu(dx). 
    \end{aligned}
    \end{equation}
    Note that
    \begin{equation*}
    \begin{aligned}
        & \int_{(0,\infty) \times (0,1)} \frac{|\log(1-\theta)|}{\lambda \theta} \, \pi_{\lambda, \theta}(d\lambda, d\theta) \\
        & \quad \leq |\log(1/2)| \int_{(0,\infty) \times (0,1/2]} \frac{1}{\lambda \theta} \, \pi_{\lambda, \theta}(d\lambda, d\theta) + 2 \int_{(0,\infty) \times (1/2,1)} \frac{|\log(1-\theta)|}{\lambda} \, \pi_{\lambda, \theta}(d\lambda, d\theta),
    \end{aligned}
    \end{equation*}
    which is finite by $\eqref{e:proofII:1}$ and \eqref{eq:typeII:cond}. Together with \eqref{e:proofII:3} this shows that all terms in \eqref{e:proofII:J22} are finite.
        
    If there is a Gaussian component $\Sigma>0$, then \eqref{eq:condII:2} holds since
    \begin{align*}
        \int_{(0,\infty) \times (0,1)} \int_{\R} f(\lambda, \theta, s)^2  \, ds \, \pi_{\lambda, \theta}(d\lambda, d\theta) &= \frac{1}{2} \int_{(0,\infty) \times (0,1)} \frac{1}{\lambda^3 \theta(1+\theta)} \, \pi_{\lambda, \theta}(d\lambda, d\theta) \\
        & \leq  \frac{1}{2} \int_{(0,\infty) \times (0,1)} \frac{1}{\lambda^3 \theta} \, \pi_{\lambda, \theta}(d\lambda, d\theta) < \infty.
    \end{align*}
    
    To verify \eqref{eq:condII:3}, we begin by decomposing the integral into two parts
    \begin{align*}
        & \int_{(0,\infty) \times (0,1)}\int_{\R} \int_{\R} \left( 1 \wedge (x^2 f(\lambda, \theta, s)^2)\right) \, \mu(dx) \, ds \, \pi_{\lambda, \theta}(d\lambda, d\theta) \\
        & \quad = \int_{(0,\infty) \times (0,1)}\int_{\R} \int_{|x| \leq 1} \left( 1 \wedge (x^2 f(\lambda, \theta, s)^2)\right) \, \mu(dx) \, ds \, \pi_{\lambda, \theta}(d\lambda, d\theta) \\
        & \quad \quad + \int_{(0,\infty) \times (0,1)}\int_{\R} \int_{|x|>1} \left( 1 \wedge (x^2 f(\lambda, \theta, s)^2)\right) \, \mu(dx) \, ds \, \pi_{\lambda, \theta}(d\lambda, d\theta) =: I_1 + I_2.
    \end{align*}
    For the first term, observe that 
    \begin{align*}
        I_1 & \leq \int_{(0,\infty) \times (0,1)}\int_{\R} \int_{|x| \leq 1} x^2 f(\lambda, \theta, s)^2 \, \mu(dx) \, ds \, \pi_{\lambda, \theta}(d\lambda, d\theta) \\
        & = \frac{1}{2} \int_{(0,\infty) \times (0,1)} \frac{1}{\lambda^3 \theta(1+\theta)} \, \pi_{\lambda, \theta}(d\lambda, d\theta) \int_{|x| \leq 1} x^2 \, \mu(dx) \\
        & \leq \frac{1}{2} \int_{(0,\infty) \times (0,1)} \frac{1}{\lambda^3 \theta} \, \pi_{\lambda, \theta}(d\lambda, d\theta) \int_{|x| \leq 1} x^2 \, \mu(dx) < \infty,
    \end{align*}
    while for the second we split the integral as
    \begin{align*}
        I_2 &= \int_{(0,\infty) \times (0,1)} \int_{\R} \int_{|x|>1} x^2 f(\lambda, \theta, s)^2 \1(|x|f(\lambda, \theta, s) \leq 1) \, \mu(dx) \, ds \, \pi_{\lambda, \theta}(d\lambda, d\theta)\\
        &\quad + \int_{(0,\infty) \times (0,1)} \int_{\R} \int_{|x|>1} \1(|x|f(\lambda, \theta, s) > 1) \, \mu(dx) \, ds \, \pi_{\lambda, \theta}(d\lambda, d\theta) =: I_{2,1} + I_{2,2}.
    \end{align*}
    Since $I_{2,1}$ can be bounded by $J_2$ defined in \eqref{e:proofII:J2}, it is finite. To establish the finiteness of $I_{2,2}$, note that the inequality $|x|f(\lambda, \theta, s) > 1$ implies $|x| \frac{e^{-\lambda \theta s}}{\lambda(1-\theta)} > 1$, which yields
    \begin{equation*}
        \1(|x|f(\lambda, \theta, s) > 1) \leq \1(|x| > \lambda(1-\theta)e^{\lambda \theta s}).
    \end{equation*}
    Consequently, the inequality can hold only for $s<s_0$, with $s_0$ given by \eqref{e:proofII:s0}. The value $s_0$ is positive whenever $|x|>\lambda(1-\theta)$, which is automatically satisfied under the above condition. Hence, as in \eqref{e:proofII:J22}, we get
    \begin{align*}
        I_{2,2} &\leq \int_{(0,\infty) \times (0,1)} \int_0^{s_0} \int_{|x|>1} \1(|x| > \lambda(1-\theta)e^{\lambda \theta s}) \, \mu(dx) \, ds \, \pi_{\lambda, \theta}(d\lambda, d\theta) \\
        &\leq \int_{(0,\infty) \times (0,1)} \int_{|x|>1} \frac{\log \left(|x| / (\lambda(1-\theta))\right)}{\lambda \theta} \1(|x| > \lambda(1-\theta)) \, \mu(dx) \, \pi_{\lambda, \theta}(d\lambda, d\theta) \\
        & \leq \int_{(0,\infty) \times (0,1)} \frac{1}{\lambda \theta} \, \pi_{\lambda, \theta}(d\lambda, d\theta) \int_{|x|> 1} \log |x| \, \mu(dx) \\
        & \quad + \int_{(0,\infty) \times (0,1)} \frac{\log(1/\lambda)}{\lambda \theta} \, \pi_{\lambda, \theta}(d\lambda, d\theta) \int_{|x|> 1} \, \mu(dx) \\
        & \quad \quad + \int_{(0,\infty) \times (0,1)} \frac{|\log(1-\theta)|}{\lambda \theta} \, \pi_{\lambda, \theta}(d\lambda, d\theta) \int_{|x|> 1} \, \mu(dx) < \infty.
    \end{align*}
    By \cite[Proposition~2.6]{rajput1989spectral}, we obtain \eqref{eq:cum:typeII}. From this identity, stationarity follows. Furthermore, from \cite[Theorem~2.7]{rajput1989spectral} it follows that $X_{\II}(t)$ is infinitely divisible with characteristic triplet $(\gamma_{\II}, \Sigma_{\II}, \mu_{\II})$.
\end{proof}

\begin{proposition}\label{prop:mc:II}
        If 
        \begin{equation*}
            \int_{|x|>1} |x|^p \, \mu(dx) < \infty, \quad \int_{(0,\infty) \times (0,1)} \frac{1}{\lambda^{p+1}\theta} \, \pi_{\lambda, \theta}(d\lambda, d\theta) < \infty,
        \end{equation*}
        then $\E |X_{\II}(0)|^p < \infty$.
\end{proposition}

\begin{proof}
    Let
    \begin{equation*}
        I_{\lambda,\theta} := \int_0^\infty \left(\frac{e^{-\lambda\theta s}-e^{-\lambda s}}{\lambda(1-\theta)}\right)^p ds.
    \end{equation*}
    Substituting $u=\lambda s$, we obtain
    \begin{equation*}
        I_{\lambda,\theta} = \frac{1}{\lambda^{p+1}(1-\theta)^p} \int_0^\infty \left(e^{-\theta u}-e^{-u}\right)^p du =: \frac{1}{\lambda^{p+1}(1-\theta)^p} J(\theta).
    \end{equation*}
    A closed form for $J(\theta)$ is given by
    \begin{equation*}
        J(\theta) = \frac{1}{1-\theta} \B\left( \frac{p \theta}{1-\theta}, p+1 \right) = \frac{1}{1-\theta} \frac{\Gamma\left( \frac{p\theta}{1-\theta}\right) \Gamma(p+1)}{\Gamma \left( \frac{p\theta}{1-\theta} +p+1 \right)},
    \end{equation*}
    where $\B$ denotes the beta function. Using the asymptotic relation $\Gamma(s) \sim 1/s$, as $s \to 0$, we obtain
    \begin{equation*}
        J(\theta) \sim \frac{1}{p \theta}, \quad \theta \to 0.
    \end{equation*}
    Similarly, since $\Gamma(s)/\Gamma(s+a) \sim s^{-a}$ as $s \to \infty$ (see e.g.~\cite{abramowitz1964handbook}),
    \begin{equation*}
        J(\theta) \sim \frac{\Gamma(p+1)}{p^{p+1}} (1-\theta)^p, \quad \theta \to 1.
    \end{equation*}
    Consequently,
    \begin{equation*}
        \lim_{\theta \to 0} \theta J(\theta) = \frac{1}{p}, \quad  \lim_{\theta \to 1} J(\theta) = 0.
    \end{equation*}
    Thus, for $\theta \in (0, 1/2]$ there exists a constant $C_0 < \infty$ such that
    \begin{equation*}
        J(\theta) \leq \frac{C_0}{\theta},
    \end{equation*}
    and therefore
    \begin{equation*}
        I_{\lambda,\theta} \leq \frac{C_0}{\lambda^{p+1} \theta (1-\theta)^p} \leq \frac{2^p C_0}{\lambda^{p+1} \theta}.
    \end{equation*}
    For $\theta \in (1/2,1)$, there exists a constant $C_1 < \infty$ such that
    \begin{equation*}
        J(\theta) \leq C_1 (1-\theta)^p.
    \end{equation*}
    Hence
    \begin{equation*}
        I_{\lambda,\theta} \leq \frac{C_1}{\lambda^{p+1}} \leq \frac{2C_1}{\lambda^{p+1} \theta}.
    \end{equation*}
    Setting
    \begin{equation*}
        C_p := \max \{ 2^p C_0, 2C_1 \},
    \end{equation*}
    we obtain the uniform bound
    \begin{equation}\label{eq:proofII:fp}
        I_{\lambda,\theta} \leq \frac{C_p}{\lambda^{p+1}\theta}, \quad \theta \in (0,1).
    \end{equation}
    For $0 < p \leq 2$, using \eqref{eq:proofII:fp}, we obtain
    \begin{align*}
        \int_{|x|>1} |x|^p \, \mu_{\II}(dx) &= \int_{(0,\infty) \times (0,1)} \int_{\R} \int_{\R} f(\lambda, \theta, s)^p |x|^p \1 (|x| f(\lambda, \theta, s)>1 ) \, \mu(dx) \, ds \, \pi_{\lambda, \theta}(d\lambda, d\theta) \\
        & \leq \int_{(0,\infty) \times (0,1)} \int_{\R} \int_{|x| \leq 1} f(\lambda, \theta, s)^2 |x|^2 \1 (|x| f(\lambda, \theta, s)>1 ) \, \mu(dx) \, ds \, \pi_{\lambda, \theta}(d\lambda, d\theta) \\
        & \quad + \int_{(0,\infty) \times (0,1)} \int_{\R} \int_{|x| > 1} f(\lambda, \theta, s)^p |x|^p \1 (|x| f(\lambda, \theta, s)>1 ) \, \mu(dx) \, ds \, \pi_{\lambda, \theta}(d\lambda, d\theta) \\
        & \leq \int_{(0,\infty) \times (0,1)} \frac{1}{\lambda^3 \theta(1+\theta)} \, \pi_{\lambda, \theta}(d\lambda, d\theta) \int_{|x| \leq 1} |x|^2 \, \mu(dx) \\
        & \quad + C_p \int_{(0,\infty) \times (0,1)} \frac{1}{\lambda^{p+1}\theta} \, \pi_{\lambda, \theta}(d\lambda, d\theta) \int_{|x| > 1} |x|^p \, \mu(dx) < \infty.
    \end{align*}
    For $p>2$, since
    \begin{equation*}
        \int_{\R} |x|^p \, \mu(dx) \leq \left( \int_{|x|\leq 1} |x|^2 \, \mu(dx) + \int_{|x|>1} |x|^p \, \mu(dx) \right) < \infty,
    \end{equation*}
    it follows that
    \begin{align*}
        \int_{|x|>1} |x|^p \, \mu_{\II}(dx) & = \int_{(0,\infty) \times (0,1)} \int_{\R} \int_{\R} f(\lambda, \theta, s)^p |x|^p \1 (|x| f(\lambda, \theta, s)>1 ) \, \mu(dx) \, ds \, \pi_{\lambda, \theta}(d\lambda, d\theta) \\
        & \leq \int_{(0,\infty) \times (0,1)} \int_{\R} f(\lambda, \theta, s)^p \, ds \, \pi_{\lambda, \theta}(d\lambda, d\theta) \int_{\R} |x|^p \, \mu(dx) \\
        & \leq C_p \int_{(0,\infty) \times (0,1)} \frac{1}{\lambda^{p+1}\theta} \, \pi_{\lambda, \theta}(d\lambda, d\theta) \int_{\R} |x|^p \, \mu(dx) < \infty.
    \end{align*}
\end{proof}

\begin{proposition}\label{prop:corr:typeII}
    If $\int_{|x|>1} x^2 \, \mu(dx) < \infty$, then the process $X_{\II}$ has finite second moment and we have
    \begin{align*}
        \E(X_{\II}(0)) &= \E(L(1)) \int_{(0,\infty) \times (0,1)} \frac{1}{\lambda^2 \theta} \, \pi_{\lambda, \theta}(d\lambda, d\theta),\\
        \Var(X_{\II}(0)) &= \frac{1}{2} \Var(L(1)) \int_{(0,\infty) \times (0,1)} \frac{1}{\lambda^3 \theta(1+\theta)} \pi_{\lambda, \theta}(d\lambda, d\theta), \\
    	r(\tau) &= \frac{\Var(L(1))}{2 \Var(X_{\II}(0))} \int_{(0,\infty) \times (0,1)} \frac{1}{\lambda^3 (1-\theta^2)} \left( \frac{e^{-\lambda \theta \tau}}{\theta} - e^{-\lambda \tau}\right) \pi_{\lambda, \theta}(d\lambda, d\theta).
    \end{align*}
\end{proposition}

\begin{proof}
    The finiteness of the second moment follows from Proposition \ref{prop:mc:II}. Following the same steps as in the proof of Proposition \ref{prop:corr:typeI}, we first obtain from \eqref{eq:cum:typeII}
        \begin{equation*}
        \kappa_{X_{\II}}'(0) = \kappa_{L}'(0) \int_{(0,\infty) \times (0,1)} \int_{-\infty}^{t_1} \frac{e^{-\lambda \theta (t_1-s)}-e^{-\lambda  (t_1-s)}}{\lambda(1-\theta) } \, ds \, \pi_{\lambda, \theta}(d\lambda, d\theta).
    \end{equation*}
    By integrating with respect to $s$, it follows that
    \begin{equation*}
        \E(X_{\II}(0))=-i\kappa_{X_{\II}}'(0)=\E(L(1)) \int_{(0,\infty) \times (0,1)} \frac{1}{\lambda^2 \theta} \, \pi_{\lambda, \theta}(d\lambda, d\theta).
    \end{equation*}
    For the covariance, we get
    \begin{equation*}
        \Cov(X_{\II}(t_1), X_{\II}(t_2)) = \frac{1}{2} \Var(L(1)) \int_{(0,\infty) \times (0,1)} \frac{1}{\lambda^3 (1-\theta^2)} \left( \frac{e^{-\lambda \theta (t_2-t_1)}}{\theta} - e^{-\lambda (t_2-t_1)}\right) \pi_{\lambda, \theta}(d\lambda, d\theta).
    \end{equation*}
    By stationarity, for $\tau>0$ we obtain
    \begin{equation*}
        \Cov(X_{\II}(0), X_{\II}(\tau)) = \frac{1}{2} \Var(L(1)) \int_{(0,\infty) \times (0,1)} \frac{1}{\lambda^3 (1-\theta^2)} \left( \frac{e^{-\lambda \theta \tau}}{\theta} - e^{-\lambda \tau}\right) \pi_{\lambda, \theta}(d\lambda, d\theta).
    \end{equation*}
    Thus, the variance is given by
    \begin{equation*}
        \Var(X_{\II}(0)) = \frac{1}{2} \Var(L(1)) \int_{(0,\infty) \times (0,1)} \frac{1}{\lambda^3 \theta(1+\theta)} \pi_{\lambda, \theta}(d\lambda, d\theta),
    \end{equation*}
    and the correlation function takes the form
    \begin{equation*}
    	r(\tau)= \frac{\Var(L(1))}{2 \Var(X_{\II}(0))} \int_{(0,\infty) \times (0,1)} \frac{1}{\lambda^3 (1-\theta^2)} \left( \frac{e^{-\lambda \theta \tau}}{\theta} - e^{-\lambda \tau}\right) \pi_{\lambda, \theta}(d\lambda, d\theta).
    \end{equation*}
\end{proof}

If the parameters $\lambda$ and $\theta$ are independent, meaning that the probability measure $\pi_{\lambda,\theta}$ factorizes as 
\begin{equation*}
    \pi_{\lambda,\theta} = \pi_{\lambda} \times \pi_{\theta},
\end{equation*}
where $\pi_{\lambda}$ is a probability measure on $(0,\infty)$ and $\pi_{\theta}$ a probability measure on $(0,1)$, then the integrability conditions in \eqref{eq:typeII:cond} reduce to
\begin{equation}\label{eq:typeII:cond2}
    \int_{(0,\infty)} \lambda^{-3} \pi_{\lambda}(d\lambda) < \infty, \quad   \int_{(0,1)} \theta^{-1} \, \pi_{\theta}(d\theta) < \infty, \quad  \int_{\left(\frac{1}{2},1\right)} |\log(1-\theta)| \, \pi_{\theta}(d\theta) < \infty.
\end{equation}
The conditions for the finiteness of moments from Proposition \ref{prop:mc:II} reduce to
\begin{equation}\label{eq:mc2:II}
    \int_{|x|>1} |x|^p \, \mu(dx) < \infty, \quad \int_{(0,\infty)} \frac{1}{\lambda^{p+1}} \, \pi_{\lambda}(d\lambda) < \infty.
\end{equation}

\begin{example}
    Suppose $X_{\II}$ is a square integrable supCAR$(2)$-II process such that $\pi_{\lambda}$ is $\Gamma(\alpha+3,1)$ distribution with $\alpha \in (0,1)$ and $\pi_{\theta}$ is $\Beta(\beta_0,\beta_1)$ distribution with $\beta_0>0$ and $\beta_1>0$, i.e.
    \begin{equation*}
        \pi_{\lambda}(d\lambda) = \frac{1}{\Gamma(\alpha+3)} \lambda^{\alpha+2} e^{-\lambda} \, d\lambda, \quad \pi_{\theta}(d\theta) = \frac{1}{\B(\beta_0, \beta_1)} \theta^{\beta_0-1} (1-\theta)^{\beta_1-1} \, d\theta.
    \end{equation*}
    The condition on $\pi_{\lambda}$ from \eqref{eq:typeII:cond2} is satisfied, since
    \begin{equation*}
        \int_{(0,\infty)} \lambda^{-3} \, \pi_{\lambda}(d\lambda) = \frac{1}{\Gamma(\alpha+3)} \int_0^{\infty} \lambda^{\alpha-1} e^{-\lambda} \, d\lambda = \frac{\Gamma(\alpha)}{\Gamma(\alpha+3)} < \infty,
    \end{equation*}
    which holds for all $\alpha>0$. Thus, under this choice of Gamma distribution, the second condition in \eqref{eq:mc2:II} is automatically satisfied, ensuring that the correlation function in Proposition \ref{prop:corr:typeII} is well-defined. It remains to determine the values of $\beta_0$ and $\beta_1$ for which $\pi_{\theta}$ satisfies the remaining conditions in \eqref{eq:typeII:cond2}. We therefore need to check the finiteness of
    \begin{equation*}
        \int_{(0,1)} \theta^{-1} \, \pi_{\theta}(d\theta) \quad \text{and} \quad \int_{\left(\frac{1}{2},1\right)} |\log(1-\theta)|\,\pi_{\theta}(d\theta).
    \end{equation*}
    For the first integral, we obtain
    \begin{equation*}
        \int_{(0,1)} \theta^{-1}\,\pi_{\theta}(d\theta) = \frac{1}{\B(\beta_0,\beta_1)} \int_0^1 \theta^{\beta_0-2} (1-\theta)^{\beta_1-1} \, d\theta = \frac{\B(\beta_0-1,\beta_1)}{\B(\beta_0,\beta_1)},
    \end{equation*}
    which requires $\beta_0>1$, since the integrand behaves as $\theta^{\beta_0-2}$ near $0$. Using this condition, the second integral can be bounded as
    \begin{align*}
        \int_{\left(\frac{1}{2},1\right)} |\log(1-\theta)|\,\pi_{\theta}(d\theta) &= \frac{1}{\B(\beta_0,\beta_1)} \int_{\frac{1}{2}}^1 |\log(1-\theta)| \theta^{\beta_0-1} (1-\theta)^{\beta_1-1} \, d\theta\\
        & \leq -\frac{1}{\B(\beta_0,\beta_1)} \int_{\frac{1}{2}}^1 \log(1-\theta) (1-\theta)^{\beta_1-1} \, d\theta\\
        &= \frac{\beta_1 \log 2 + 1}{2^{\beta_1} \beta_1^2 \B(\beta_0,\beta_1)} < \infty,
    \end{align*}
    so the second condition in \eqref{eq:typeII:cond2} is satisfied for all $\beta_1>0$.

    By Fubini's theorem and integration with respect to $\lambda$, we obtain
    \begin{align*}
        r(\tau)&= \frac{\Var(L(1))}{2 \Gamma(\alpha+3) \Var(X_{\II}(0))} \int_{(0,1)} \int_{0}^{\infty} \frac{1}{\lambda^3 (1-\theta^2)} \left( \frac{e^{-\lambda \theta \tau}}{\theta} - e^{-\lambda \tau}\right)\lambda^{\alpha+2}e^{-\lambda}\,d\lambda\,\pi_{\theta}(d\theta)\\
        &= \frac{\Gamma(\alpha) \Var(L(1))}{2 \Gamma(\alpha+3) \Var(X_{\II}(0))} \int_{(0,1)} \frac{1}{1-\theta^2} \left( \frac{1}{\theta(\theta \tau+1)^{\alpha}} - \frac{1}{(\tau+1)^{\alpha}}\right) \,\pi_{\theta}(d\theta) \\
        &= \frac{\Var(L(1))}{2 \B(\beta_0,\beta_1) \Gamma(\alpha+3) \Var(X_{\II}(0))} \int_0^1 \frac{\theta^{\beta_0-1} (1-\theta)^{\beta_1-2}}{1+\theta} \left( \frac{1}{\theta(\theta \tau+1)^{\alpha}} - \frac{1}{(\tau+1)^{\alpha}}\right) d\theta.
    \end{align*}
    
    To analyze the asymptotic behavior of $r(\tau)$ as $\tau \to \infty$, we additionally assume that $\beta_0 > \alpha + 1$. Define
    \begin{equation*}
        g(\tau,\theta) = \frac{\theta^{\beta_0-1} (1-\theta)^{\beta_1-2}}{1+\theta} \left( \frac{1}{\theta(\theta \tau+1)^{\alpha}} - \frac{1}{(\tau+1)^{\alpha}}\right), \quad \theta \in (0,1), \ \tau>0.
    \end{equation*}
    For each fixed $\theta \in (0,1)$,
    \begin{equation*}
        (\theta\tau+1)^{-\alpha} \sim \tau^{-\alpha} \theta^{-\alpha}, \quad (\tau+1)^{-\alpha} \sim \tau^{-\alpha}, \quad \text{as }\tau \to \infty,
    \end{equation*}
    and therefore
    \begin{equation}\label{eq:ex2:asymp}
        g(\tau,\theta) \sim \tau^{-\alpha} \frac{\theta^{\beta_{0}-1}(1-\theta)^{\beta_{1}-2}}{1+\theta} \left(\theta^{-(\alpha+1)}-1\right) =: \tau^{-\alpha} h(\theta), \quad \text{as }\tau \to \infty.
    \end{equation}
    Hence $\tau\mapsto g(\tau,\theta)$ is regularly varying of index $-\alpha$. Let $f(\tau,\theta) := \tau^{\alpha} g(\tau,\theta).$ By \eqref{eq:ex2:asymp}, for each fixed $\theta \in (0,1)$,
    \begin{equation*}
        f(\tau,\theta) \to h(\theta), \quad \text{as }\tau \to \infty.
    \end{equation*}
    For $\theta \in (0, 1/2]$ one easily checks that
    \begin{equation*}
        |f(\tau,\theta)| \leq \left(\theta^{\beta_0-\alpha-2}+ \theta^{\beta_0-1}\right)(1-\theta)^{\beta_1-2} \leq \theta^{\beta_0-\alpha-2}+ \theta^{\beta_0-1}.
    \end{equation*}
    For fixed $\tau$, the function
    \begin{equation*}
        k_{\tau}(\theta) := \frac{\tau^{\alpha}}{\theta(\theta \tau+1)^{\alpha}}
    \end{equation*}
    is continuously differentiable on $(1/2,1)$. By the mean value theorem, for each $\theta \in (1/2,1)$ there exists $\xi \in (\theta,1)$ such that
    \begin{equation*}
        k_{\tau}(\theta) - k_{\tau}(1) = k_{\tau}'(\xi)(\theta-1).
    \end{equation*}
    Hence
    \begin{equation*}
        f(\tau,\theta) = -\frac{\theta^{\beta_0-1}(1-\theta)^{\beta_1-1}}{1+\theta} k_{\tau}'(\xi).
    \end{equation*}
    Since $\xi \in (1/2,1]$, all factors in $k_{\tau}'(\xi)$ are bounded and there exists $C>0$ such that $|k_{\tau}'(\xi)| \leq C$. Consequently,
    \begin{equation*}
        |f(\tau,\theta)| \leq C, \quad \theta \in (1/2,1).
    \end{equation*}
    Since
    \begin{equation*}
        \int_{(0,1)} \theta^{\beta_0-\alpha-2} \, d\theta < \infty,
    \end{equation*}
    whenever $\beta_0 > \alpha+1$, the dominated convergence theorem yields
    \begin{equation*}
        \int_0^1 f(\tau, \theta) \, d\theta \to \int_0^1 h(\theta) \, d\theta, \quad \tau \to \infty.
    \end{equation*}
    
    By \cite[Theorem~4.1.4]{bingham1989regular}, it follows that
    \begin{align*}
        \tau^{\alpha} r(\tau) &= \frac{\Var(L(1))}{2\B(\beta_0,\beta_1)\Gamma(\alpha+3)\Var(X_{\II}(0))} \int_0^1 f(\tau,\theta)\,d\theta \\
        & \quad \sim \frac{\Var(L(1))}{2\B(\beta_0,\beta_1)\Gamma(\alpha+3)\Var(X_{\II}(0))} \int_0^1 h(\theta)\,d\theta.
    \end{align*}
    Hence,
    \begin{equation*}
        r(\tau) \sim C(\alpha,\beta_0,\beta_1) \tau^{-\alpha}, \quad \tau\to\infty,
    \end{equation*}
    where
    \begin{equation*}
        C(\alpha,\beta_0,\beta_1) = \frac{\Var(L(1))}{2\B(\beta_0,\beta_1)\Gamma(\alpha+3)\Var(X_{\II}(0))} \int_0^1 \frac{\theta^{\beta_{0}-1}(1-\theta)^{\beta_{1}-2}}{1+\theta} \left(\theta^{-(\alpha+1)}-1\right)\,d\theta.
    \end{equation*}
    In particular, whenever $\alpha\in (0,1]$, the correlation function of the supCAR(2)-II process exhibits long-range dependence.
\end{example}

\section{SupCAR$(2)$ of type III}\label{sec:typeIII}

Finally, we turn to the case when the eigenvalues \eqref{e:eigenvalues} are complex conjugate and the kernel function is $g_{\III}$ from \eqref{e:gfunctions}. Similarly as for the supCAR$(2)$-II process, instead of modeling the distribution of elements of matrix $A$, it is more convenient to specify the distribution of its eigenvalues. The complex conjugate eigenvalues can always be written as
\begin{equation*}
    \xi_1 = r e^{i\psi}, \quad \xi_2 = r e^{-i\psi},
\end{equation*}
for $r>0$ and $\psi \in (\frac{\pi}{2},\pi)$, since the eigenvalues should have negative real part. Therefore we introduce the following definition which introduces randomness through parameters $r$ and $\psi$.

\begin{definition}
Let $\pi_{r,\psi}$ be a probability measure on $(0,\infty) \times (\frac{\pi}{2}, \pi)$ and $\Lambda_{\III}$ a L\'evy basis on $(0,\infty) \times (\frac{\pi}{2}, \pi) \times \R$ corresponding to the L\'evy-Khintchine triplet $(\gamma, \Sigma, \mu)$ and control measure $\pi_{r,\psi} \times \Leb$. Assume that $\int_{|x|>1} \log |x| \, \mu(dx)$ and
\begin{equation}\label{eq:typeIII:cond}
\begin{aligned}
    &\int_{(0,\infty) \times (\frac{\pi}{2},\pi)} \frac{1}{r^3 |\cos \psi|} \, \pi_{r, \psi}(dr, d\psi) < \infty, \\
    &\int_{(0,\infty) \times (\frac{\pi}{2},\pi)} \frac{\log (\sin \psi)}{r \cos \psi} \, \pi_{r, \psi}(dr, d\psi) < \infty, \\
    &\int_{(\frac{\pi}{2},\pi)} \frac{1}{|\cos \psi|} \, \pi_{\psi}(d\psi) < \infty.
\end{aligned}
\end{equation}
A \textbf{supCAR$(2)$ process of type III} (supCAR$(2)$-III) is defined by
\begin{equation*}
X_{\III}(t) = \int_{(0,\infty)\times (\frac{\pi}{2},\pi) \times \R} \frac{e^{r(t-s)\cos \psi}}{r \sin \psi} \sin \left( r(t-s)\sin \psi \right) \1_{[0,\infty)}(t-s)\, \Lambda_{\III}(dr, d\psi, ds).
\end{equation*}
\end{definition}

\begin{theorem}
supCAR(2)-III process is well-defined and stationary with cumulant function of its finite-dimensional distributions given by
    \begin{align}\label{eq:cum:typeIII}
        &\kappa_{\left(X_{\III}(t_1), \dots, X_{\III}(t_m)\right)}(\zeta_1, \dots \zeta_m) \notag \\
        &\quad = \int_{(0,\infty) \times (\frac{\pi}{2},\pi)} \int_{\R} \kappa_{L} \left( \sum_{j=1}^m \zeta_j \frac{e^{r(t_j-s)\cos \psi}}{r \sin \psi} \sin \left( r(t_j-s)\sin \psi \right) \1_{[0,\infty)}(t_j-s) \right) ds \, \pi_{r, \psi}(dr, d\psi),
    \end{align}
where $t_1 < \cdots < t_m$. The distribution of $X_{\III}(t)$ is infinitely divisible with characteristic triplet $(\gamma_{\III}, \Sigma_{\III}, \mu_{\III})$ given by
\begin{align*}
    \gamma_{\III} &= \gamma \int_{(0,\infty) \times (\frac{\pi}{2},\pi)} \frac{1}{r^2} \, \pi_{r, \psi}(dr, d\psi) \\
    & \quad + \int_{(0,\infty) \times (\frac{\pi}{2},\pi)} \int_{0}^{\infty} \int_{\R} \Big[ \tau\left (x\frac{e^{rs \cos \psi}}{r \sin \psi} \sin \left( rs \sin \psi \right)\right)\\
    &\quad \quad \quad \quad \quad \quad \quad \quad \quad \quad \quad \quad \quad \quad -\frac{e^{rs \cos \psi}}{r \sin \psi} \sin \left( rs \sin \psi \right)\tau(x) \Big] \, \mu(dx) \, ds \, \pi_{r, \psi}(dr, d\psi),\\
    \Sigma_{\III} &= \frac{\Sigma}{4} \int_{(0,\infty) \times (\frac{\pi}{2},\pi)} \frac{1}{r^3|\cos \psi|} \, \pi_{r, \psi}(dr, d\psi), \\
    \mu_{\III}(B) &= \int_{(0,\infty) \times (\frac{\pi}{2},\pi)} \int_0^{\infty} \int_{\R} \1_{B}\left(\frac{e^{rs \cos \psi}}{r \sin \psi} \sin \left( rs \sin \psi \right) x\right) \, \mu(dx) \, ds \, \pi_{r, \psi}(dr, d\psi), \quad B \in \mathcal{B}(\R).
\end{align*}
\end{theorem}

\begin{proof}
    Let
    \begin{equation*}
        f(r, \psi, s) := \frac{e^{rs \cos \psi}}{r \sin \psi} \sin \left( rs \sin \psi \right) \1_{[0,\infty)}(s).
    \end{equation*}
    We now verify the conditions \eqref{eq:rajros:mma}, which in the present setting take the form
    \begin{align}
        & \int_{(0,\infty) \times (\frac{\pi}{2},\pi)} \int_{\R} \Big| \gamma f(r, \psi, s) + \int_{\R} \left( \tau(xf(r, \psi, s))-f(r, \psi, s)\tau(x) \right) \, \mu(dx) \Big| \, ds \, \pi_{r, \psi}(dr, d\psi) < \infty, \label{eq:condIII:1}\\
        & \int_{(0,\infty) \times (\frac{\pi}{2},\pi)} \int_{\R} f(r, \psi, s)^2  \, ds \, \pi_{r, \psi}(dr, d\psi) < \infty, \label{eq:condIII:2}\\
        & \int_{(0,\infty) \times (\frac{\pi}{2},\pi)}\int_{\R} \int_{\R} \left( 1 \wedge (x^2 f(r, \psi, s)^2)\right) \, \mu(dx) \, ds \, \pi_{r, \psi}(dr, d\psi) < \infty. \label{eq:condIII:3}
    \end{align}
    First, note that for $x>0$, we have $\coth(x) \leq 1 + 1/x$, and therefore
    \begin{equation*}
        \coth \left(\frac{\pi |\cot \psi |}{2}\right) \leq 1 + \frac{2}{\pi | \cos \psi |}.    
    \end{equation*}
    Hence
    \begin{align*}
         & \int_{(0,\infty) \times (\frac{\pi}{2},\pi)} \frac{1}{r^2} \coth \left(\frac{\pi |\cot \psi |}{2}\right) \, \pi_{r, \psi}(dr, d\psi) \\
         & \quad \leq  \int_{(0,\infty) \times (\frac{\pi}{2},\pi)} \frac{1}{r^2} \, \pi_{r, \psi}(dr, d\psi) + \frac{2}{\pi} \int_{(0,\infty) \times (\frac{\pi}{2},\pi)} \frac{1}{r^2 |\cos \psi|} \, \pi_{r, \psi}(dr, d\psi).
    \end{align*}
    Since $r^{-2} \leq r^{-3}$ for $r\in(0,1)$ and $r^{-2} \leq 1$ for $r\geq1$, and because $\pi_{r, \psi}$ is a probability measure,
    \begin{align*}
        \int_{(0,\infty) \times (\frac{\pi}{2},\pi)} \frac{1}{r^2} \, \pi_{r, \psi}(dr, d\psi) & = \int_{(0,1) \times (\frac{\pi}{2},\pi)} \frac{1}{r^2} \, \pi_{r, \psi}(dr, d\psi) + \int_{[1,\infty) \times (\frac{\pi}{2},\pi)} \frac{1}{r^2} \, \pi_{r, \psi}(dr, d\psi) \\
        & \leq \int_{(0,1) \times (\frac{\pi}{2},\pi)} \frac{1}{r^3} \, \pi_{r, \psi}(dr, d\psi) + \int_{[1,\infty) \times (\frac{\pi}{2},\pi)} \, \pi_{r, \psi}(dr, d\psi) \\
        & \leq \int_{(0,1) \times (\frac{\pi}{2},\pi)} \frac{1}{r^3 |\cos \psi|} \, \pi_{r, \psi}(dr, d\psi) + 1 < \infty
    \end{align*}
    and
    \begin{equation}\label{eq:proofIII:r2cos}
        \begin{aligned}
            &\int_{(0,\infty) \times (\frac{\pi}{2},\pi)} \frac{1}{r^2 |\cos \psi|} \, \pi_{r, \psi}(dr, d\psi) \\
            & \quad = \int_{(0,1) \times (\frac{\pi}{2},\pi)} \frac{1}{r^2 |\cos \psi|} \, \pi_{r, \psi}(dr, d\psi) + \int_{[1,\infty) \times (\frac{\pi}{2},\pi)} \frac{1}{r^2 |\cos \psi|} \, \pi_{r, \psi}(dr, d\psi)\\
            & \quad \leq \int_{(0,1) \times (\frac{\pi}{2},\pi)} \frac{1}{r^3 |\cos \psi|} \, \pi_{r, \psi}(dr, d\psi) + \int_{(\frac{\pi}{2},\pi)} \frac{1}{|\cos \psi|} \, \pi_{\psi}(d\psi) < \infty.
        \end{aligned}
        \end{equation}
    Therefore
    \begin{equation}\label{eq:proofIII:coth}
        \int_{(0,\infty) \times (\frac{\pi}{2},\pi)} \frac{1}{r^2} \coth \left(\frac{\pi |\cot \psi |}{2}\right) \, \pi_{r, \psi}(dr, d\psi) < \infty.
    \end{equation}
    Note that
    \begin{equation}\label{eq:proofIII:rcos}
        \int_{(0,\infty) \times (\frac{\pi}{2},\pi)} \frac{1}{r |\cos \psi|} \, \pi_{r, \psi}(dr, d\psi) < \infty,
    \end{equation}
    by the same argument as in \eqref{eq:proofIII:r2cos}, with $r^{-1}$ estimated by $r^{-3}$ on $(0,1)$ and by $1$ on $[1,\infty)$. Similarly, using $|\log r|/r \leq r^{-3}$ for $r\in(0,1)$ and $|\log r|/r \leq 1$ for $r\geq1$, we can estimate
    \begin{equation}\label{eq:proofIII:logr}
    \begin{aligned}
        &\int_{(0,\infty) \times (\frac{\pi}{2},\pi)} \frac{\log r}{r\cos \psi} \, \pi_{r, \psi}(dr, d\psi) \\
        &\quad \leq \int_{(0,\infty) \times (\frac{\pi}{2},\pi)} \frac{|\log r|}{r|\cos \psi|} \, \pi_{r, \psi}(dr, d\psi) \\
        &\quad \leq \int_{(0,1) \times (\frac{\pi}{2},\pi)} \frac{|\log r|}{r|\cos \psi|} \, \pi_{r, \psi}(dr, d\psi) + \int_{[1,\infty) \times (\frac{\pi}{2},\pi)} \frac{|\log r|}{r|\cos \psi|} \, \pi_{r, \psi}(dr, d\psi) \\
        &\quad \leq \int_{(0,1) \times (\frac{\pi}{2},\pi)} \frac{1}{r^3|\cos \psi|} \, \pi_{r, \psi}(dr, d\psi) + \int_{(\frac{\pi}{2},\pi)} \frac{1}{r|\cos \psi|} \, \pi_{\psi}(d\psi) < \infty.
    \end{aligned}
    \end{equation}
    
    To verify \eqref{eq:condIII:1}, we first consider the integral
    \begin{align*} 
        & \int_{(0,\infty) \times (\frac{\pi}{2},\pi)} \int_{\R} |f(r, \psi, s)| \, ds \, \pi_{r, \psi}(dr, d\psi) \\
        & \quad = \int_{(0,\infty) \times (\frac{\pi}{2},\pi)} \int_0^{\infty} \frac{e^{rs \cos \psi}}{r \sin \psi} \left| \sin \left( rs \sin \psi \right) \right| \, ds \, \pi_{r, \psi}(dr, d\psi) \\
        & \quad = \int_{(0,\infty) \times (\frac{\pi}{2},\pi)} \frac{1}{r^2} \coth \left(\frac{\pi |\cot \psi |}{2}\right) \, \pi_{r, \psi}(dr, d\psi),
    \end{align*}
    which is finite by \eqref{eq:proofIII:coth}. For the second term in \eqref{eq:condIII:1}, we have
    \begin{align*}
        &\int_{(0,\infty) \times (\frac{\pi}{2},\pi)} \int_{\R} \int_{\R} \left| \tau \left(xf(r, \psi, s)\right)-f(r, \psi, s)\tau(x) \right| \, \mu(dx) \, ds \, \pi_{r, \psi}(dr, d\psi)\\
        &\quad = \int_{(0,\infty) \times (\frac{\pi}{2},\pi)} \int_{\R} \int_{|x|\leq 1} |xf(r, \psi, s)| \1(|xf(r, \psi, s)|> 1) \, \mu(dx) \, ds \, \pi_{r, \psi}(dr, d\psi)\\
        &\quad \quad + \int_{(0,\infty) \times (\frac{\pi}{2},\pi)} \int_{\R} \int_{|x|> 1} |xf(r, \psi, s)| \1(|xf(r, \psi, s)| \leq 1)\, \mu(dx) \, ds \, \pi_{r, \psi}(dr, d\psi)\\
        &\quad =: J_1+J_2.
    \end{align*}
    For $J_1$, it follows that
    \begin{align*}
        J_1 & \leq  \int_{(0,\infty) \times (\frac{\pi}{2},\pi)} \int_{\R} \int_{|x|\leq 1} x^2 f(r, \psi, s)^2 \1(|xf(r, \psi, s)|> 1) \, \mu(dx) \, ds \, \pi_{r, \psi}(dr, d\psi)\\
        & \leq \frac{1}{4} \int_{(0,\infty) \times (\frac{\pi}{2},\pi)} \frac{1}{r^3 |\cos \psi|} \, \pi_{r, \psi}(dr, d\psi) \int_{|x|\leq 1} x^2 \, \mu(dx) < \infty.
    \end{align*}
    Since $|\sin(rs\sin\psi)| \leq 1$ for all $s>0$, the inequality 
    \begin{equation}\label{eq:ineq:III}
        |x|\frac{e^{rs \cos \psi}}{r \sin \psi} |\sin(rs\sin\psi)| \leq 1
    \end{equation}
    is satisfied whenever $|x|\frac{e^{rs \cos \psi}}{r \sin \psi} \leq 1$. Define
    \begin{equation*}
        s_0 := \frac{\log(r\sin \psi / |x|)}{r \cos \psi}.
    \end{equation*}
    Then, for $s>s_0$, the inequality \eqref{eq:ineq:III} will certainly hold. If $|x| \leq r \sin \psi$, then $s<s_0$ and \eqref{eq:ineq:III} is satisfied for every $s>0$. Thus, we decompose $J_2$ as
    \begin{align*}
        J_2 &= \int_{(0,\infty) \times (\frac{\pi}{2},\pi)} \int_{\R} \int_{|x|> 1} |xf(r, \psi, s)| \1(|xf(r, \psi, s)| \leq 1) \1(|x| \leq r \sin \psi) \, \mu(dx) \, ds \, \pi_{r, \psi}(dr, d\psi)\\
        & \quad + \int_{(0,\infty) \times (\frac{\pi}{2},\pi)} \int_{\R} \int_{|x|> 1} |xf(r, \psi, s)| \1(|xf(r, \psi, s)| \leq 1) \1(|x| > r \sin \psi) \, \mu(dx) \, ds \, \pi_{r, \psi}(dr, d\psi) \\
        & =: J_{2,1} + J_{2,2}.
    \end{align*}
    For $J_{2,1}$, we have
    \begin{align*}
        J_{2,1} &= \int_{(0,\infty) \times (\frac{\pi}{2},\pi)} \int_{\R} \int_{|x|> 1} |xf(r, \psi, s)| \1(|x| \leq r \sin \psi) \, \mu(dx) \, ds \, \pi_{r, \psi}(dr, d\psi) \\
        & \leq \int_{(0,\infty) \times (\frac{\pi}{2},\pi)} \frac{1}{r^2 |\sin 2\psi|} \int_{|x|> 1} |x| \1(|x| \leq r \sin \psi) \, \mu(dx) \, \pi_{r, \psi}(dr, d\psi) \\
        & \leq \int_{(0,\infty) \times (\frac{\pi}{2},\pi)} \frac{1}{2 r |\cos \psi|} \1(r \sin \psi >1) \, \pi_{r, \psi}(dr, d\psi) \int_{|x|> 1} |x| \, \mu(dx) < \infty,
    \end{align*}
    and finiteness follows from \eqref{eq:proofIII:rcos}. For $J_{2,2}$ we obtain
    \begin{align*}
        J_{2,2} &\leq \int_{(0,\infty) \times (\frac{\pi}{2},\pi)} \int_0^{s_0} \int_{|x|> 1} \1(|x| > r \sin \psi) \, \mu(dx) \, ds \, \pi_{r, \psi}(dr, d\psi)\\
        & \quad + \int_{(0,\infty) \times (\frac{\pi}{2},\pi)} \int_{s_0}^{\infty} \int_{|x|> 1} |xf(r, \psi, s)| \, \mu(dx) \, ds \, \pi_{r, \psi}(dr, d\psi) \\
        & \leq \int_{(0,\infty) \times (\frac{\pi}{2},\pi)} \int_{|x|> 1} \frac{\log(r\sin \psi / |x|)}{r \cos \psi} \1(|x| > r \sin \psi) \, \mu(dx) \, \pi_{r, \psi}(dr, d\psi)\\
        & \quad + \int_{(0,\infty) \times (\frac{\pi}{2},\pi)} \int_{|x|> 1} \frac{1}{r |\cos \psi|} \, \mu(dx) \, \pi_{r, \psi}(dr, d\psi) \\
        & \leq \int_{(0,\infty) \times (\frac{\pi}{2},\pi)} \frac{|\log r|}{r |\cos \psi|} \, \pi_{r, \psi}(dr, d\psi) \int_{|x|> 1} \, \mu(dx) \\
        & \quad + \int_{(0,\infty) \times (\frac{\pi}{2},\pi)} \frac{\log(\sin \psi)}{r \cos \psi} \, \pi_{r, \psi}(dr, d\psi) \int_{|x|> 1} \, \mu(dx) \\
        & \quad \quad + \int_{(0,\infty) \times (\frac{\pi}{2},\pi)} \frac{1}{r |\cos \psi|} \, \pi_{r, \psi}(dr, d\psi) \int_{|x|> 1} \log |x| \, \mu(dx) \\
        & \quad \quad \quad \int_{(0,\infty) \times (\frac{\pi}{2},\pi)} \frac{1}{r |\cos \psi|} \, \pi_{r, \psi}(dr, d\psi) \int_{|x|> 1} \, \mu(dx).
    \end{align*}
    The finiteness of the above expression follows from \eqref{eq:typeIII:cond}, \eqref{eq:proofIII:rcos} and \eqref{eq:proofIII:logr}.
    
    If there is a Gaussian component $\Sigma>0$, then \eqref{eq:condIII:2} holds since
    \begin{equation*}
        \int_{(0,\infty) \times (\frac{\pi}{2},\pi)} \int_{\R} f(r, \psi, s)^2  \, ds \, \pi_{r, \psi}(dr, d\psi) = \frac{1}{4} \int_{(0,\infty) \times (\frac{\pi}{2},\pi)} \frac{1}{r^3|\cos \psi|} \, \pi_{r, \psi}(dr, d\psi) < \infty.
    \end{equation*}
    
    In order to show \eqref{eq:condIII:3}, we split the integral as follows
    \begin{align*}
        & \int_{(0,\infty) \times (\frac{\pi}{2},\pi)}\int_{\R} \int_{\R} \left( 1 \wedge (x^2 f(r, \psi, s)^2)\right) \, \mu(dx) \, ds \, \pi_{r, \psi}(dr, d\psi) \\
        & \quad = \int_{(0,\infty) \times (\frac{\pi}{2},\pi)}\int_{\R} \int_{|x| \leq 1} \left( 1 \wedge (x^2 f(r, \psi, s)^2)\right) \, \mu(dx) \, ds \, \pi_{r, \psi}(dr, d\psi) \\
        & \quad \quad + \int_{(0,\infty) \times (\frac{\pi}{2},\pi)}\int_{\R} \int_{|x|>1} \left( 1 \wedge (x^2 f(r, \psi, s)^2)\right) \, \mu(dx) \, ds \, \pi_{r, \psi}(dr, d\psi) =: I_1 + I_2.
    \end{align*}
    For $I_1$, we get 
    \begin{align*}
        I_1 & \leq \int_{(0,\infty) \times (\frac{\pi}{2},\pi)}\int_{\R} \int_{|x| \leq 1} x^2 f(r, \psi, s)^2 \, \mu(dx) \, ds \, \pi_{r, \psi}(dr, d\psi) \\
        & = \frac{1}{4} \int_{(0,\infty) \times (\frac{\pi}{2},\pi)} \frac{1}{r^3 |\cos \psi|} \, \pi_{r, \psi}(dr, d\psi) \int_{|x| \leq 1} x^2 \, \mu(dx) < \infty,
    \end{align*}
    while for $I_2$ we split the integral as
    \begin{align*}
        I_2 &= \int_{(0,\infty) \times (\frac{\pi}{2},\pi)} \int_{\R} \int_{|x|>1} x^2 f(r, \psi, s)^2 \1(|x f(r, \psi, s)| \leq 1) \, \mu(dx) \, ds \, \pi_{r, \psi}(dr, d\psi)\\
        &\quad + \int_{(0,\infty) \times (\frac{\pi}{2},\pi)} \int_{\R} \int_{|x|>1} \1(|x f(r, \psi, s)| > 1) \, \mu(dx) \, ds \, \pi_{r, \psi}(dr, d\psi) =: I_{2,1} + I_{2,2}.
    \end{align*}
    Note that $I_{2,1} \leq J_2 < \infty$. To show that $I_{2,2}$ is finite, first note that the inequality $|x f(r, \psi, s)| > 1$ implies $|x| \frac{e^{rs \cos \psi}}{r \sin \psi} > 1$, i.e.~it holds that
    \begin{equation*}
        \1(|x f(r, \psi, s)| > 1) \leq \1(|x| > r \sin \psi e^{-rs \cos \psi}).
    \end{equation*}
    In that case, it must hold that $s<s_0$, so we obtain
    \begin{align*}
        I_{2,2} &\leq \int_{(0,\infty) \times (\frac{\pi}{2},\pi)} \int_0^{s_0} \int_{|x|>1} \1(|x| > r \sin \psi e^{-rs \cos \psi}) \, \mu(dx) \, ds \, \pi_{r, \psi}(dr, d\psi)\\
        & \leq \int_{(0,\infty) \times (\frac{\pi}{2},\pi)} \int_{|x|>1} \frac{\log(r \sin \psi / |x|)}{r \cos \psi} \, \mu(dx)  \,  \pi_{r, \psi}(dr, d\psi) \\
        & = \int_{(0,\infty) \times (\frac{\pi}{2},\pi)} \frac{\log r}{r\cos \psi} \, \pi_{r, \psi}(dr, d\psi) \int_{|x|> 1} \, \mu(dx)\\
        & \quad + \int_{(0,\infty) \times (\frac{\pi}{2},\pi)} \frac{\log(\sin \psi)}{r\cos \psi} \, \pi_{r, \psi}(dr, d\psi) \int_{|x|> 1} \, \mu(dx) \\
        & \quad \quad + \int_{(0,\infty) \times (\frac{\pi}{2},\pi)} \frac{1}{r|\cos \psi|} \, \pi_{r, \psi}(dr, d\psi) \int_{|x|> 1} \log |x| \, \mu(dx) < \infty.
    \end{align*}
    As in the previous theorems, \eqref{eq:cum:typeIII} together with the representations of $(\gamma_{\III}, \Sigma_{\III}, \mu_{\III})$ follow from \cite[Proposition~2.6]{rajput1989spectral} and \cite[Theorem~2.7]{rajput1989spectral}.
\end{proof}

\begin{proposition}\label{prop:mc:III}
    If
    \begin{equation*}
        \int_{|x|>1} |x|^p \, \mu(dx) < \infty, \quad \int_{(0,\infty) \times (\frac{\pi}{2},\pi)} \frac{1}{r^{p+1} |\cos \psi|} \, \pi_{r, \psi}(dr, d\psi) < \infty,
    \end{equation*}
    then $\E |X_{\III}(0)|^p < \infty$.
\end{proposition}

\begin{proof}
    Define
    \begin{equation*}
        I_{r,\psi} := \int_0^{\infty} \left( \frac{e^{r s \cos \psi}}{r \sin \psi} |\sin(r s \sin \psi)| \right)^p ds.
    \end{equation*}
    Substituting $u = r s \sin \psi$, we obtain
    \begin{equation*}
        I_{r,\psi} = \frac{1}{r^{p+1} \sin^{p+1}\psi} \int_0^{\infty} e^{p u \cot \psi} |\sin u|^p \, du =: \frac{1}{r^{p+1} \sin^{p+1}\psi} J(\psi).
    \end{equation*}
    Since $p \cot \psi < 0$ for $\psi \in (\frac{\pi}{2}, \pi)$, we consider two cases: $p \cot \psi \in [-1,0)$ and $p \cot \psi \in (-\infty, -1)$. If $p \cot \psi \in [-1,0)$, we have
    \begin{equation*}
        J(\psi) \leq \int_0^{\infty} e^{p u \cot \psi} \, du = \frac{\sin\psi}{p |\cos\psi|},
    \end{equation*}
    which implies
    \begin{equation*}
        I_{r,\psi} \leq \frac{1}{p r^{p+1} \sin^p \psi |\cos\psi|}.
    \end{equation*}
    Since $|\cot \psi| \leq 1/p$, there exists a positive constant $C_0$ such that $\sin \psi \geq C_0$. Therefore, we obtain
    \begin{equation*}
        I_{r,\psi} \leq \frac{1}{p C_0^p r^{p+1} |\cos\psi|}.
    \end{equation*}
    Now we consider the second case, i.e.~$p \cot \psi \in (-\infty,-1)$. Using the inequality $|\sin u| \leq |u|$ for all $u \in \R$, we get
    \begin{equation*}
        J(\psi) \leq \int_0^{\infty} e^{p u \cot \psi} |u|^p \, du = \frac{\Gamma(p+1)}{p^{p+1} |\cot \psi|^{p+1}}.
    \end{equation*}
    Hence,
    \begin{equation*}
        I_{r,\psi} \leq \frac{\Gamma(p+1)}{p^{p+1} r^{p+1} \sin^{p+1} \psi |\cot \psi|^{p+1}} = \frac{\Gamma(p+1)}{p^{p+1} r^{p+1} |\cos \psi|^{p+1}}.
    \end{equation*}
    Since $|\cot \psi| > 1/p$, there exists a positive constant $C_1$ such that $|\cos \psi| \geq C_1$, and thus
    \begin{equation*}
        I_{r,\psi} \leq \frac{\Gamma(p+1)}{p^{p+1}C_1^p r^{p+1} |\cos \psi|}.
    \end{equation*}
    Setting
    \begin{equation*}
        C_p := \max \left \{ \frac{1}{p C_0^p}, \frac{\Gamma(p+1)}{p^{p+1} C_1^p}\right \},
    \end{equation*}
    we obtain the uniform bound
    \begin{equation}\label{eq:proofIII:fp}
        I_{r,\psi} \leq \frac{C_p}{r^{p+1} |\cos\psi|}, \quad \psi \in (\pi/2, \pi).
    \end{equation}
    For $0 < p \leq 2$, we have
    \begin{align*}
        \int_{|x|>1} |x|^p \, \mu_{\III}(dx) &= \int_{(0,\infty) \times (\frac{\pi}{2},\pi)} \int_{\R} \int_{\R} |x f(r, \psi, s)|^p \1 (|x f(r, \psi, s)|>1 ) \, \mu(dx) \, ds \, \pi_{r, \psi}(dr, d\psi) \\
        & \leq \int_{(0,\infty) \times (\frac{\pi}{2},\pi)} \int_{\R} \int_{|x| \leq 1} f(r, \psi, s)^2 |x|^2 \1 (|x f(r, \psi, s)|>1 ) \, \mu(dx) \, ds \, \pi_{r, \psi}(dr, d\psi) \\
        & \quad + \int_{(0,\infty) \times (\frac{\pi}{2},\pi)} \int_{\R} \int_{|x| > 1} |f(r, \psi, s)|^p |x|^p \1 (|x f(r, \psi, s)|>1 ) \, \mu(dx) \, ds \, \pi_{r, \psi}(dr, d\psi) \\
        & \leq \frac{1}{4} \int_{(0,\infty) \times (\frac{\pi}{2},\pi)} \frac{1}{r^3 |\cos \psi|} \, \pi_{r, \psi}(dr, d\psi) \int_{|x| \leq 1} |x|^2 \, \mu(dx) \\
        & \quad + C_p \int_{(0,\infty) \times (\frac{\pi}{2},\pi)} \frac{1}{r^{p+1} |\cos \psi|} \, \pi_{r, \psi}(dr, d\psi) \int_{|x| > 1} |x|^p \, \mu(dx) < \infty.
    \end{align*}
    Using
    \begin{equation*}
        \int_{\R} |x|^p \, \mu(dx) \leq\int_{|x|\leq 1} |x|^2 \, \mu(dx) + \int_{|x|>1} |x|^p \, \mu(dx) < \infty,
    \end{equation*}
    we obtain
    \begin{align*}
        \int_{|x|>1} |x|^p \, \mu_{\II}(dx) &= \int_{(0,\infty) \times (\frac{\pi}{2},\pi)} \int_{\R} \int_{\R} |f(r, \psi, s)|^p |x|^p \1 (|x f(r, \psi, s)|>1 ) \, \mu(dx) \, ds \, \pi_{r, \psi}(dr, d\psi) \\
        & \leq C_p \int_{(0,\infty) \times (\frac{\pi}{2},\pi)} \frac{1}{r^{p+1} |\cos \psi|} \, \pi_{r, \psi}(dr, d\psi) \int_{\R} |x|^p \, \mu(dx) < \infty.
    \end{align*}
\end{proof}

\begin{proposition}\label{prop:corr:typeIII}
    If $\int_{|x|>1} x^2 \, \mu(dx) < \infty$, then the process $X_{\III}$ has finite second moment and we have
    \begin{align*}
        \E(X_{\III}(0)) &= \E(L(1)) \int_{(0,\infty) \times (\frac{\pi}{2},\pi)} \frac{1}{r^2} \, \pi_{r, \psi}(dr, d\psi),\\
        \Var(X_{\III}(0)) &= \frac{1}{4} \Var(L(1)) \int_{(0,\infty) \times (\frac{\pi}{2},\pi)} \frac{1}{r^3 |\cos \psi|} \, \pi_{r, \psi}(dr, d\psi), \\
    	r(\tau) &= \frac{\Var(L(1))}{2\Var(X_{\III}(0))} \int_{(0,\infty) \times (\frac{\pi}{2},\pi)} \frac{e^{r\tau \cos \psi}}{r^3 \sin 2\psi} \sin\left( r\tau \sin \psi - \psi\right) \pi_{r, \psi}(dr, d\psi).
    \end{align*}
\end{proposition}

\begin{proof}
    The finiteness of the second moment follows from Proposition \ref{prop:mc:III}. In a similar way as for the previous two types, from \eqref{eq:cum:typeIII} we get
    \begin{equation*}
        \kappa_{X_{\III}}'(0) = \kappa_{L}'(0) \int_{(0,\infty) \times (\frac{\pi}{2},\pi)} \int_{-\infty}^{t_1} \frac{e^{r(t_1-s)\cos \psi}}{r \sin \psi} \sin \left( r(t_1-s)\sin \psi \right) ds \, \pi_{r, \psi}(dr, d\psi).
    \end{equation*}
    Integrating with respect to $s$ yields
    \begin{equation*}
        \E(X_{\III}(0))=-i\kappa_{X_{\III}}'(0)= \E(L(1)) \int_{(0,\infty) \times (\frac{\pi}{2},\pi)} \frac{1}{r^2} \, \pi_{r, \psi}(dr, d\psi).
    \end{equation*}
    For the autocovariance function, we obtain
    \begin{equation*}
        \Cov(X_{\III}(t_1),X_{\III}(t_2))=\frac{1}{2} \Var(L(1)) \int_{(0,\infty) \times (\frac{\pi}{2},\pi)} \frac{e^{r(t_2-t_1) \cos \psi}}{r^3 \sin 2\psi} \sin\left( r(t_2-t_1) \sin \psi - \psi\right) \pi_{r, \psi}(dr, d\psi).
    \end{equation*}
    Hence
    \begin{equation*}
        \Cov(X_{\III}(0),X_{\III}(\tau))=\frac{1}{2} \Var(L(1)) \int_{(0,\infty) \times (\frac{\pi}{2},\pi)} \frac{e^{r\tau \cos \psi}}{r^3 \sin 2\psi} \sin\left( r\tau \sin \psi - \psi\right) \pi_{r, \psi}(dr, d\psi).
    \end{equation*}
    Then, the variance is
    \begin{equation*}
        \Var(X_{\III}(0))=\frac{1}{4} \Var(L(1)) \int_{(0,\infty) \times (\frac{\pi}{2},\pi)} \frac{1}{r^3 |\cos \psi|} \, \pi_{r, \psi}(dr, d\psi)
    \end{equation*}
    and the correlation function is 
    \begin{equation*}
        r(\tau)=\frac{\Var(L(1))}{2\Var(X_{\III}(0))} \int_{(0,\infty) \times (\frac{\pi}{2},\pi)} \frac{e^{r\tau \cos \psi}}{r^3 \sin 2\psi} \sin\left( r\tau \sin \psi - \psi\right) \pi_{r, \psi}(dr, d\psi).
    \end{equation*}
\end{proof}

If the parameters $r$ and $\psi$ are independent, so that the probability measure $\pi_{r,\psi}$ factorizes as
\begin{equation*}
    \pi_{r,\psi} = \pi_{r} \times \pi_{\psi},
\end{equation*}
with $\pi_{r}$ a probability measure on $(0,\infty)$ and $\pi_{\psi}$ a
probability measure on $(\frac{\pi}{2},\pi)$, then the integrability
conditions in \eqref{eq:typeIII:cond} take the form
\begin{equation}\label{eq:typeIII:cond2}
    \int_{(0,\infty)} r^{-3} \,\pi_{r}(dr) < \infty, \quad \int_{(\frac{\pi}{2},\pi)} \frac{1}{|\cos \psi|} \, \pi_{\psi}(d\psi) < \infty, \quad \int_{(\frac{\pi}{2},\pi)} \frac{\log(\sin \psi)}{\cos \psi} \, \pi_{\psi}(d\psi) < \infty.
\end{equation}
The conditions for the finiteness of moments from Proposition \ref{prop:mc:III} reduce to
\begin{equation*}
    \int_{|x|>1} |x|^p \, \mu(dx) < \infty, \quad \int_{(0,\infty)} r^{-(p+1)} \, \pi_{r}(dr) < \infty.
\end{equation*}

\begin{figure}[!ht]
    \centering

    \begin{subfigure}{0.48\textwidth}
        \centering
        \includegraphics[width=\linewidth]{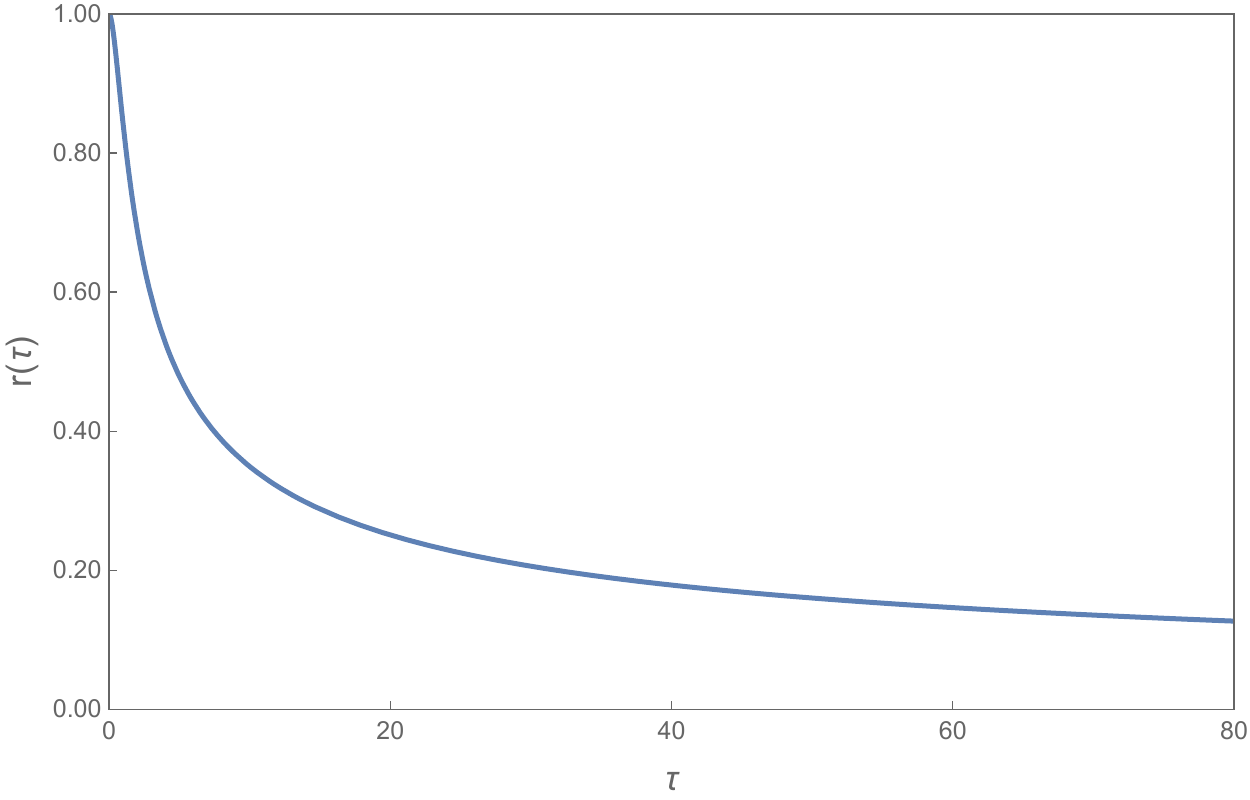}
        \caption{$\pi_r \sim \Gamma(3.5,1)$}
    \end{subfigure}\hfill
    \begin{subfigure}{0.48\textwidth}
        \centering
        \includegraphics[width=\linewidth]{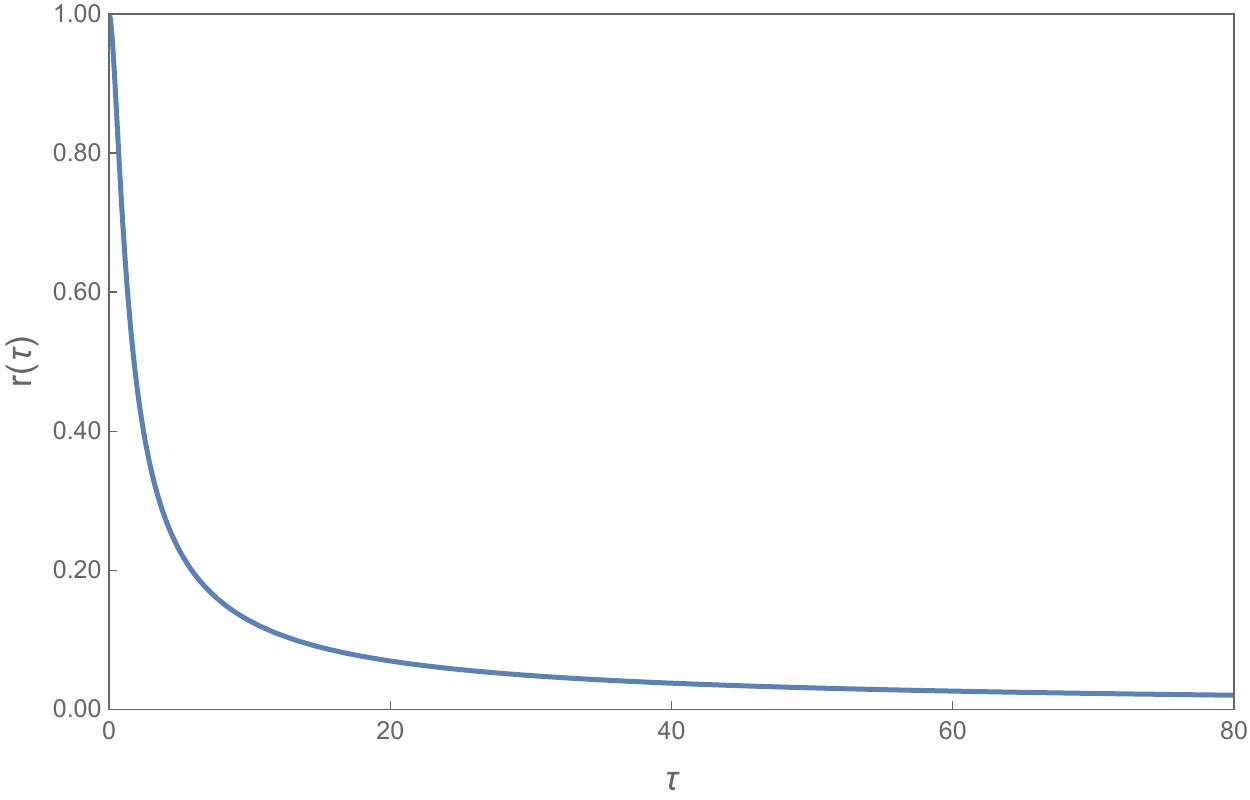}
        \caption{$\pi_r \sim \Gamma(3.9,1)$}
    \end{subfigure}

    \vspace{0.5cm}
    
    \begin{subfigure}{0.48\textwidth}
        \centering
        \includegraphics[width=\linewidth]{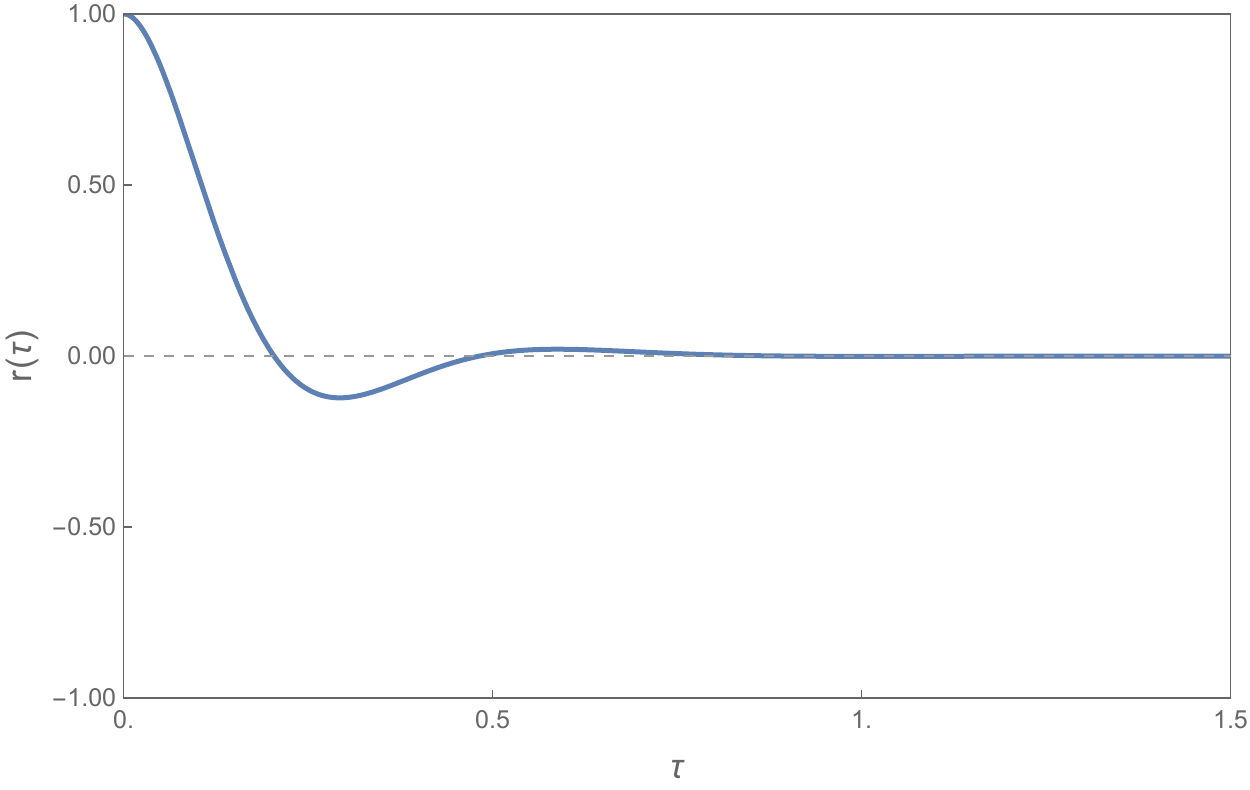}
        \caption{$\pi_r \sim \Gamma(15,1)$}
    \end{subfigure}\hfill
    \begin{subfigure}{0.48\textwidth}
        \centering
        \includegraphics[width=\linewidth]{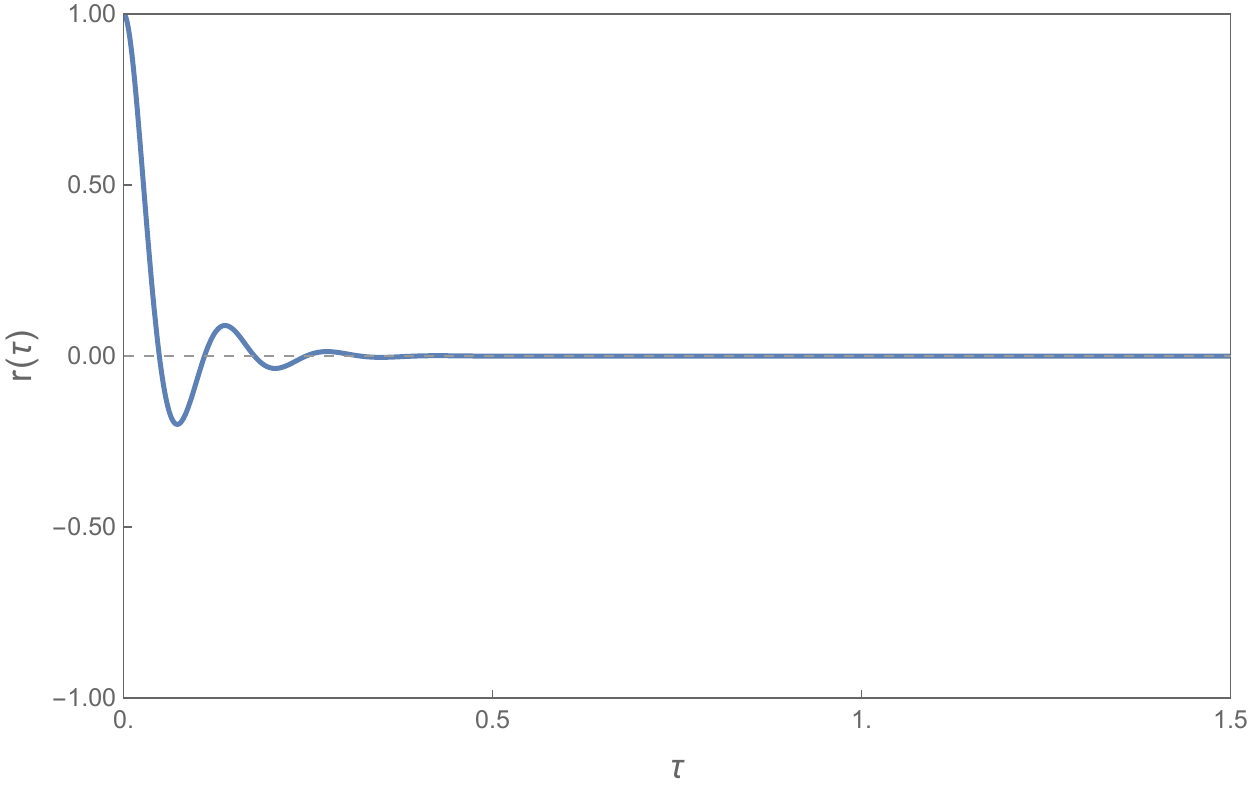}
        \caption{$\pi_r \sim \Gamma(50,1)$}
    \end{subfigure}

    \caption{Correlation functions of supCAR$(2)$-III processes for different choices of $\pi_r$ and with $\pi_{\psi}$ given by \eqref{eq:exIII:pipsi}. For shape parameters $\alpha+3 \in (3,4]$, the correlation function exhibits non-oscillatory long-range dependence, while for larger values of the shape parameter oscillations arise.}
    \label{fig:III}
\end{figure}

\begin{example}
    Suppose $X_{\III}$ is a square integrable supCAR$(2)$-III process such that $\pi_{r}$ is $\Gamma(\alpha+3,1)$ distribution with $\alpha > 0$ and $\pi_{\psi}$ is given by
    \begin{equation}\label{eq:exIII:pipsi}
        \pi_{\psi}(d\psi) = - \sin 2\psi \, \1_{(\frac{\pi}{2},\pi)}(\psi).
    \end{equation}
    The condition on $\pi_{r}$ in \eqref{eq:typeIII:cond2} holds, since
    \begin{equation*}
        \int_{(0,\infty)} r^{-3} \,\pi_{r}(dr) = \frac{1}{\Gamma(\alpha+3)} \int_0^{\infty} r^{\alpha-1}e^{-r} \, dr = \frac{\Gamma(\alpha)}{\Gamma(\alpha+3)} < \infty.
    \end{equation*}
    The remaining conditions in \eqref{eq:typeIII:cond2} are satisfied as follows
    \begin{equation*}
        \int_{(\frac{\pi}{2},\pi)} \frac{1}{|\cos \psi|} \, \pi_{\psi}(d\psi) = 2 \int_{\frac{\pi}{2}}^{\pi} \sin \psi \, d\psi = 2 < \infty,
    \end{equation*}
    and
    \begin{equation*}
        \int_{(\frac{\pi}{2},\pi)} \frac{\log(\sin \psi)}{\cos \psi} \, \pi_{\psi}(d\psi) = -2\int_{\frac{\pi}{2}}^{\pi} \log(\sin \psi) \sin \psi \, d\psi = 2-\log 4 < \infty.
    \end{equation*}
    In particular, $X_{\III}$ has finite variance and the correlation function is well-defined.
    
    We now analyze the asymptotic behavior of $r(\tau)$ for large lags $\tau$. Applying the Fubini–Tonelli theorem, we obtain
    \begin{align*}
        r(\tau) &= \frac{\Gamma(\alpha+3)}{\Gamma(\alpha)} \int_{(\frac{\pi}{2},\pi)} \int_{(0,\infty)} \frac{e^{r\tau \cos \psi}}{r^3 \sin 2\psi} \sin\left( r\tau \sin \psi - \psi\right) \pi_{r}(dr) \, \pi_{\psi}(d\psi) \\
        & = -\frac{1}{\Gamma(\alpha)} \int_{\frac{\pi}{2}}^{\pi} \int_0^{\infty} e^{r\tau \cos \psi} \sin\left( r\tau \sin \psi - \psi\right) r^{\alpha-1} e^{-r} \, dr \, d\psi.
    \end{align*}
    Note that
    \begin{equation*}
        e^{r\tau \cos \psi} \sin\left( r\tau \sin \psi - \psi\right) = \Im (e^{-i \psi} e^{r \tau e^{i \psi}} ),
    \end{equation*}
    and
    \begin{equation*}
        \int_0^{\infty} t^{z-1} e ^{-at} \, dt = a^{-z} \Gamma(z).
    \end{equation*}
    Hence
    \begin{align*}
        r(\tau) &= -\frac{1}{\Gamma(\alpha)} \int_{\frac{\pi}{2}}^{\pi} \Im \left( e^{-i \psi} \int_0^{\infty} r^{\alpha-1} e^{-r(1-\tau e^{i \psi})} \, dr \right) \, d\psi  \\
        &= -\int_{\frac{\pi}{2}}^{\pi} \Im \left( e^{-i \psi} (1-\tau e^{i \psi})^{-\alpha} \right) \, d\psi.
    \end{align*}
    Define
    \begin{equation*}
        g(\tau, \psi) := \Im \left( e^{-i \psi} (1-\tau e^{i \psi})^{-\alpha} \right), \quad \psi \in \left(\frac{\pi}{2},\pi \right), \ \tau > 0.
    \end{equation*}
    To analyze the asymptotic behavior of $1-\tau e^{i \psi}$, write 
    \begin{equation*}
        1-\tau e^{i \psi} = \rho e^{i \varphi}, \quad \rho>0, \ \varphi \in [-\pi, \pi).
    \end{equation*}
    Since
    \begin{equation*}
        1-\tau e^{i \psi} = 1 - \tau \cos \psi - i \tau \sin \psi,
    \end{equation*}
    we have
    \begin{equation*}
        \rho = \sqrt{1 - 2 \tau \cos \psi + \tau^2}, \quad \varphi = \arctan \left( \frac{-\tau \sin \psi}{1-\tau \cos \psi} \right),
    \end{equation*}
    and
    \begin{equation*}
        \Re(1-\tau e^{i \psi}) = 1 - \tau \cos \psi > 0, \quad \Im(1-\tau e^{i \psi}) = - \tau \sin \psi < 0.
    \end{equation*}
    Thus $1-\tau e^{i\psi}$ eventually lies in the fourth quadrant. As $\tau \to \infty$, we have
    \begin{equation*}
        \rho \sim \tau, \quad \varphi \to \psi - \pi.
    \end{equation*}
    Hence
    \begin{equation*}
        e^{-i \psi} (1-\tau e^{i \psi})^{-\alpha} \sim \tau^{-\alpha} e^{i(\alpha(\pi-\psi) - \psi)}, \quad \tau \to \infty,
    \end{equation*}
    and therefore
    \begin{equation}\label{eq:ex4:asymp}
        g(\tau, \psi) \sim \tau^{-\alpha} \sin\left( \alpha(\pi-\psi)-\psi \right) =: \tau^{-\alpha} h(\psi), \quad \tau \to \infty.
    \end{equation}
    Thus $\tau \mapsto g(\tau, \psi)$ is regularly varying of index $-\alpha$. Let $f(\tau,\psi) := \tau^{\alpha} g(\tau, \psi)$. By \eqref{eq:ex4:asymp}, $f(\tau,\psi) \to h(\psi)$ as $\tau\to\infty$, and
    \begin{equation*}
        |f(\tau,\psi)| \leq 1.
    \end{equation*}
    Hence, by dominated convergence,
    \begin{equation*}
        \int_{\frac{\pi}{2}}^{\pi} f(\tau, \psi) \, d\psi \to \int_{\frac{\pi}{2}}^{\pi} h(\psi) \, d\psi, \quad \tau \to \infty.
    \end{equation*}
    Observe that
    \begin{equation*}
        \int_{\frac{\pi}{2}}^{\pi} h(\psi) \, d\psi = - \frac{1+\sin \frac{\alpha \pi}{2}}{1+\alpha}.
    \end{equation*}
    Finally, \cite[Theorem~4.1.4]{bingham1989regular} yields 
    \begin{equation*}
        r(\tau) \sim \frac{1+\sin \frac{\alpha \pi}{2}}{1+\alpha} \tau^{-\alpha}, \quad \tau \to \infty.
    \end{equation*}
    In particular, the supCAR(2)-III process exhibits long-range dependence for $\alpha\in(0,1]$. Moreover, it can exhibit oscillatory and non-monotonic behavior as shown in Figure \ref{fig:III}, where correlation functions are plotted for different choices of $\pi_r$.
\end{example}

\bibliographystyle{apamc}
\bibliography{References}

@article {rajput1989spectral,
    AUTHOR = {Rajput, Balram S. and Rosi\'{n}ski, Jan},
     TITLE = {Spectral representations of infinitely divisible processes},
   JOURNAL = {Probab. Theory Related Fields},
  FJOURNAL = {Probability Theory and Related Fields},
    VOLUME = {82},
      YEAR = {1989},
    NUMBER = {3},
     PAGES = {451--487},
      ISSN = {0178-8051,1432-2064},
   MRCLASS = {60H05 (60E07 60G57 60J30)},
  MRNUMBER = {1001524},
MRREVIEWER = {John\ P.\ Nolan},
       DOI = {10.1007/BF00339998},
       URL = {https://doi.org/10.1007/BF00339998},
}

@article{bn2001,
author = {Barndorff-Nielsen, O. E.},
title = {Superposition of {O}rnstein--{U}hlenbeck type processes},
journal = {Theory of Probability \& Its Applications},
volume = {45},
number = {2},
pages = {175-194},
year = {2001},
doi = {10.1137/S0040585X97978166},
URL = {https://doi.org/10.1137/S0040585X97978166},
eprint = {https://doi.org/10.1137/S0040585X97978166}
}

@book {barndorff2018ambit,
    AUTHOR = {Barndorff-Nielsen, Ole E. and Benth, Fred Espen and Veraart,
              Almut E. D.},
     TITLE = {Ambit stochastics},
    SERIES = {Probability Theory and Stochastic Modelling},
    VOLUME = {88},
 PUBLISHER = {Springer, Cham},
      YEAR = {2018},
     PAGES = {xxv+402},
      ISBN = {978-3-319-94128-8; 978-3-319-94129-5},
   MRCLASS = {60-02 (60G60 60Gxx 60Hxx 62Mxx 65Cxx 91G20)},
  MRNUMBER = {3839270},
MRREVIEWER = {Anthony\ R\'{e}veillac},
       DOI = {10.1007/978-3-319-94129-5},
       URL = {https://doi.org/10.1007/978-3-319-94129-5},
}

@article {barndorff2011multivariate,
    AUTHOR = {Barndorff-Nielsen, Ole Eiler and Stelzer, Robert},
     TITLE = {Multivariate sup{OU} processes},
   JOURNAL = {Ann. Appl. Probab.},
  FJOURNAL = {The Annals of Applied Probability},
    VOLUME = {21},
      YEAR = {2011},
    NUMBER = {1},
     PAGES = {140--182},
      ISSN = {1050-5164,2168-8737},
   MRCLASS = {60J75 (60E07 60G51 60G57)},
  MRNUMBER = {2759198},
MRREVIEWER = {Yasushi\ Ishikawa},
       DOI = {10.1214/10-AAP690},
       URL = {https://doi.org/10.1214/10-AAP690},
}

@article {barndorff2013multivariate,
    AUTHOR = {Barndorff-Nielsen, Ole Eiler and Stelzer, Robert},
     TITLE = {The multivariate sup{OU} stochastic volatility model},
   JOURNAL = {Math. Finance},
  FJOURNAL = {Mathematical Finance. An International Journal of Mathematics,
              Statistics and Financial Economics},
    VOLUME = {23},
      YEAR = {2013},
    NUMBER = {2},
     PAGES = {275--296},
      ISSN = {0960-1627,1467-9965},
   MRCLASS = {91B70 (60G51 60J70)},
  MRNUMBER = {3034078},
MRREVIEWER = {George\ Stoica},
       DOI = {10.1111/j.1467-9965.2011.00494.x},
       URL = {https://doi.org/10.1111/j.1467-9965.2011.00494.x},
}

@book {bingham1989regular,
    AUTHOR = {Bingham, N. H. and Goldie, C. M. and Teugels, J. L.},
     TITLE = {Regular variation},
    SERIES = {Encyclopedia of {M}athematics and its {A}pplications},
    VOLUME = {27},
 PUBLISHER = {Cambridge University Press, Cambridge},
      YEAR = {1989},
     PAGES = {xx+494},
      ISBN = {0-521-37943-1},
   MRCLASS = {26A12 (11-01)},
  MRNUMBER = {1015093},
}

@book {abramowitz1964handbook,
    AUTHOR = {Abramowitz, Milton and Stegun, Irene A.},
     TITLE = {Handbook of mathematical functions with formulas, graphs, and
              mathematical tables},
    SERIES = {National Bureau of Standards Applied Mathematics Series, No.
              55},
      NOTE = {For sale by the Superintendent of Documents},
 PUBLISHER = {U. S. Government Printing Office, Washington, DC},
      YEAR = {1964},
     PAGES = {xiv+1046},
   MRCLASS = {33.00 (65.05)},
  MRNUMBER = {167642},
MRREVIEWER = {D.\ H.\ Lehmer},
}

@article {brockwell2005levy,
    AUTHOR = {Brockwell, Peter and Marquardt, Tina},
     TITLE = {L\'{e}vy-driven and fractionally integrated {ARMA} processes
              with continuous time parameter},
   JOURNAL = {Statist. Sinica},
  FJOURNAL = {Statistica Sinica},
    VOLUME = {15},
      YEAR = {2005},
    NUMBER = {2},
     PAGES = {477--494},
      ISSN = {1017-0405,1996-8507},
   MRCLASS = {62M10 (60G51 62P05)},
  MRNUMBER = {2190215},
}

@incollection {brockwell2001249,
    AUTHOR = {Brockwell, P. J.},
     TITLE = {Continuous-time {ARMA} processes},
 BOOKTITLE = {Stochastic processes: theory and methods},
    SERIES = {Handbook of Statist.},
    VOLUME = {19},
     PAGES = {249--276},
 PUBLISHER = {North-Holland, Amsterdam},
      YEAR = {2001},
      ISBN = {0-444-50014-6},
   MRCLASS = {62M10},
  MRNUMBER = {1861726},
       DOI = {10.1016/S0169-7161(01)19011-5},
       URL = {https://doi.org/10.1016/S0169-7161(01)19011-5},
}

@article {brockwell2001levy,
    AUTHOR = {Brockwell, P. J.},
     TITLE = {L\'{e}vy-driven {CARMA} processes},
      NOTE = {Nonlinear non-Gaussian models and related filtering methods
              (Tokyo, 2000)},
   JOURNAL = {Ann. Inst. Statist. Math.},
  FJOURNAL = {Annals of the Institute of Statistical Mathematics},
    VOLUME = {53},
      YEAR = {2001},
    NUMBER = {1},
     PAGES = {113--124},
      ISSN = {0020-3157,1572-9052},
   MRCLASS = {62M05 (62M10)},
  MRNUMBER = {1820952},
       DOI = {10.1023/A:1017972605872},
       URL = {https://doi.org/10.1023/A:1017972605872},
}

@book {Sato,
    AUTHOR = {Sato, {Ken-iti}},
     TITLE = {L\'{e}vy processes and infinitely divisible distributions},
    SERIES = {Cambridge Studies in Advanced Mathematics},
    VOLUME = {68},
 PUBLISHER = {Cambridge University Press, Cambridge},
      YEAR = {1999},
     PAGES = {xii+486},
      ISBN = {0-521-55302-4},
   MRCLASS = {60G51 (60E07 60G18 60G52 60J45)},
  MRNUMBER = {1739520},
MRREVIEWER = {N.\ H.\ Bingham},
}

@article{marquardt2007generating,
  title={Generating long memory models based on {CARMA} processes},
  author={Marquardt, TINA and James, LANCELOT F},
  journal={Technical report},
  year={2007}
}

@article {brockwelllindner2009,
    AUTHOR = {Brockwell, Peter J. and Lindner, Alexander},
     TITLE = {Existence and uniqueness of stationary {L}\'{e}vy-driven
              {CARMA} processes},
   JOURNAL = {Stochastic Process. Appl.},
  FJOURNAL = {Stochastic Processes and their Applications},
    VOLUME = {119},
      YEAR = {2009},
    NUMBER = {8},
     PAGES = {2660--2681},
      ISSN = {0304-4149,1879-209X},
   MRCLASS = {60G51 (60G10 62M10)},
  MRNUMBER = {2532218},
MRREVIEWER = {Robert\ Stelzer},
       DOI = {10.1016/j.spa.2009.01.006},
       URL = {https://doi.org/10.1016/j.spa.2009.01.006},
}

@article {kaagstrom1977bounds,
    AUTHOR = {K\'{a}gstr\"{o}m, Bo},
     TITLE = {Bounds and perturbation bounds for the matrix exponential},
   JOURNAL = {Nordisk Tidskr. Informationsbehandling (BIT)},
  FJOURNAL = {Nordisk Tidskrift for Informationsbehandling},
    VOLUME = {17},
      YEAR = {1977},
    NUMBER = {1},
     PAGES = {39--57},
      ISSN = {0901-246X},
   MRCLASS = {65F99 (15A15)},
  MRNUMBER = {440896},
MRREVIEWER = {Rolf\ Jeltsch},
       DOI = {10.1007/bf01932398},
       URL = {https://doi.org/10.1007/bf01932398},
}

@book {samorodnitsky2016book,
    AUTHOR = {Samorodnitsky, Gennady},
     TITLE = {Stochastic processes and long range dependence},
    SERIES = {Springer Series in Operations Research and Financial
              Engineering},
 PUBLISHER = {Springer, Cham},
      YEAR = {2016},
     PAGES = {xi+415},
      ISBN = {978-3-319-45574-7; 978-3-319-45575-4},
   MRCLASS = {60G10 (60E07 60F17 60G18 60G22 60G52)},
  MRNUMBER = {3561100},
MRREVIEWER = {Yizao Wang},
       DOI = {10.1007/978-3-319-45575-4},
       URL = {https://doi.org/10.1007/978-3-319-45575-4},
}

@article{chatzigeorgiou2013bounds,
author = {Chatzigeorgiou, Ioannis},
year = {2013},
month = {08},
pages = {1505-1508},
title = {Bounds on the {L}ambert function and their application to the outage analysis of user cooperation},
volume = {17},
journal = {Communications Letters, IEEE},
doi = {10.1109/LCOMM.2013.070113.130972}
}

@book {mezo2022,
    AUTHOR = {Mez\H{o}, Istv\'{a}n},
     TITLE = {The {L}ambert {W} function---its generalizations and
              applications},
    SERIES = {Discrete Mathematics and its Applications (Boca Raton)},
 PUBLISHER = {CRC Press, Boca Raton, FL},
      YEAR = {2022},
     PAGES = {xxi+252},
      ISBN = {978-0-367-76683-2; 978-1-03-222339-1; 978-1-00-316810-2},
   MRCLASS = {33E20},
  MRNUMBER = {4600791},
}

@article{wolpert2005,
title = {Fractional {O}rnstein–{U}hlenbeck {L}\'evy processes and the {T}elecom process: {U}pstairs and downstairs},
journal = {Signal Processing},
volume = {85},
number = {8},
pages = {1523-1545},
year = {2005},
issn = {0165-1684},
doi = {https://doi.org/10.1016/j.sigpro.2004.09.016},
url = {https://www.sciencedirect.com/science/article/pii/S0165168405000757},
author = {Robert L. Wolpert and Murad S. Taqqu}
}

@incollection {bngamma2016,
    AUTHOR = {Barndorff-Nielsen, Ole E.},
     TITLE = {Gamma kernels and {BSS}/{LSS} processes},
 BOOKTITLE = {Advanced modelling in mathematical finance},
    SERIES = {Springer Proc. Math. Stat.},
    VOLUME = {189},
     PAGES = {41--61},
 PUBLISHER = {Springer, Cham},
      YEAR = {2016},
   MRCLASS = {60G07 (60G10 60G22 60G51 60H05 62M09 91G70)},
  MRNUMBER = {3630532},
       DOI = {10.1007/978-3-319-45875-5\_2}
}

@article {granger1980,
    AUTHOR = {Granger, C. W. J.},
     TITLE = {Long memory relationships and the aggregation of dynamic
              models},
   JOURNAL = {J. Econometrics},
  FJOURNAL = {Journal of Econometrics},
    VOLUME = {14},
      YEAR = {1980},
    NUMBER = {2},
     PAGES = {227--238},
      ISSN = {0304-4076},
   MRCLASS = {62M10 (62P20)},
  MRNUMBER = {597259},
MRREVIEWER = {Bert M. Steece},
       DOI = {10.1016/0304-4076(80)90092-5},
       URL = {https://doi.org/10.1016/0304-4076(80)90092-5},
}

@article {robinson1978statistical,
    AUTHOR = {Robinson, P. M.},
     TITLE = {Statistical inference for a random coefficient autoregressive
              model},
   JOURNAL = {Scand. J. Statist.},
  FJOURNAL = {Scandinavian Journal of Statistics. Theory and Applications},
    VOLUME = {5},
      YEAR = {1978},
    NUMBER = {3},
     PAGES = {163--168},
      ISSN = {0303-6898},
   MRCLASS = {62M10},
  MRNUMBER = {509453},
MRREVIEWER = {Ra\'{u}l Pedro Mentz},
}

@article {oppenheim2004,
    AUTHOR = {Oppenheim, Georges and Viano, Marie-Claude},
     TITLE = {Aggregation of random parameters {O}rnstein-{U}hlenbeck or
              {AR} processes: some convergence results},
   JOURNAL = {J. Time Ser. Anal.},
  FJOURNAL = {Journal of Time Series Analysis},
    VOLUME = {25},
      YEAR = {2004},
    NUMBER = {3},
     PAGES = {335--350},
      ISSN = {0143-9782},
   MRCLASS = {62M10 (60G10)},
  MRNUMBER = {2062677},
MRREVIEWER = {Ji\v{r}\'{\i} And\v{e}l},
       DOI = {10.1111/j.1467-9892.2004.01775.x},
       URL = {https://doi.org/10.1111/j.1467-9892.2004.01775.x},
}

@article {chong2001,
    AUTHOR = {Chong, Terence Tai-leung and Wong, Kwan-to},
     TITLE = {Time series properties of aggregated {${\rm AR}(2)$}
              processes},
   JOURNAL = {Econom. Lett.},
  FJOURNAL = {Economics Letters},
    VOLUME = {73},
      YEAR = {2001},
    NUMBER = {3},
     PAGES = {325--332},
      ISSN = {0165-1765},
   MRCLASS = {62M10},
  MRNUMBER = {1866754},
       DOI = {10.1016/S0165-1765(01)00504-3},
       URL = {https://doi.org/10.1016/S0165-1765(01)00504-3},
}

@article {philippe2014,
    AUTHOR = {Philippe, Anne and Puplinskaite, Donata and Surgailis,
              Donatas},
     TITLE = {Contemporaneous aggregation of triangular array of
              random-coefficient {AR}(1) processes},
   JOURNAL = {J. Time Series Anal.},
  FJOURNAL = {Journal of Time Series Analysis},
    VOLUME = {35},
      YEAR = {2014},
    NUMBER = {1},
     PAGES = {16--39},
      ISSN = {0143-9782},
   MRCLASS = {62M10 (62G20)},
  MRNUMBER = {3148246},
MRREVIEWER = {N. Leonenko},
       DOI = {10.1111/jtsa.12045},
       URL = {https://doi.org/10.1111/jtsa.12045},
}

@article {linden1999,
    AUTHOR = {Linden, Mikael},
     TITLE = {Time series properties of aggregated {${\rm AR}(1)$} processes
              with uniformly distributed coefficients},
   JOURNAL = {Econom. Lett.},
  FJOURNAL = {Economics Letters},
    VOLUME = {64},
      YEAR = {1999},
    NUMBER = {1},
     PAGES = {31--36},
      ISSN = {0165-1765},
   MRCLASS = {62P20 (62M10)},
  MRNUMBER = {1703637},
       DOI = {10.1016/S0165-1765(99)00072-5},
       URL = {https://doi.org/10.1016/S0165-1765(99)00072-5},
}

@article {robinson1986cross,
    AUTHOR = {Robinson, P. M.},
     TITLE = {On a model for a time series of cross-sections},
      NOTE = {Essays in time series and allied processes},
   JOURNAL = {J. Appl. Probab.},
  FJOURNAL = {Journal of Applied Probability},
      YEAR = {1986},
    VOLUME = {Special Vol. 23{\rm A}},
     PAGES = {113--125},
      ISSN = {0021-9002},
   MRCLASS = {62M10},
  MRNUMBER = {803167},
MRREVIEWER = {Tarmo Pukkila},
       DOI = {10.1017/s0021900200117024},
       URL = {https://doi.org/10.1017/s0021900200117024},
}

@article {barndorff2001nongaussian,
    AUTHOR = {Barndorff-Nielsen, Ole E. and Shephard, Neil},
     TITLE = {Non-{G}aussian {O}rnstein-{U}hlenbeck-based models and some of
              their uses in financial economics},
   JOURNAL = {J. R. Stat. Soc. Ser. B Stat. Methodol.},
  FJOURNAL = {Journal of the Royal Statistical Society. Series B.
              Statistical Methodology},
    VOLUME = {63},
      YEAR = {2001},
    NUMBER = {2},
     PAGES = {167--241},
      ISSN = {1369-7412},
   MRCLASS = {62M07 (62M09 62M10 62P20)},
  MRNUMBER = {1841412},
       DOI = {10.1111/1467-9868.00282},
       URL = {https://doi.org/10.1111/1467-9868.00282},
}

@article {moser2011tail,
    AUTHOR = {Moser, Martin and Stelzer, Robert},
     TITLE = {Tail behavior of multivariate {L}\'{e}vy-driven mixed moving
              average processes and sup{OU} stochastic volatility models},
   JOURNAL = {Adv. in Appl. Probab.},
  FJOURNAL = {Advances in Applied Probability},
    VOLUME = {43},
      YEAR = {2011},
    NUMBER = {4},
     PAGES = {1109--1135},
      ISSN = {0001-8678},
   MRCLASS = {60G51 (60G10 60G70 60H30 91B70 91G20)},
  MRNUMBER = {2867948},
MRREVIEWER = {Josep Vives},
       DOI = {10.1239/aap/1324045701},
       URL = {https://doi.org/10.1239/aap/1324045701},
}

@article {stelzer2015moment,
    AUTHOR = {Stelzer, Robert and Tosstorff, Thomas and Wittlinger, Marc},
     TITLE = {Moment based estimation of sup{OU} processes and a related
              stochastic volatility model},
   JOURNAL = {Stat. Risk Model.},
  FJOURNAL = {Statistics \& Risk Modeling with Applications in Finance and
              Insurance},
    VOLUME = {32},
      YEAR = {2015},
    NUMBER = {1},
     PAGES = {1--24},
      ISSN = {2193-1402},
   MRCLASS = {62M09 (60G51 62M10 91B25 91B84)},
  MRNUMBER = {3331779},
       DOI = {10.1515/strm-2012-1152},
       URL = {https://doi.org/10.1515/strm-2012-1152},
}

@article {puplinskaite2009,
    AUTHOR = {Puplinskait\.{e}, D. and Surgailis, D.},
     TITLE = {Aggregation of random-coefficient {AR}(1) process with
              infinite variance and common innovations},
   JOURNAL = {Lith. Math. J.},
  FJOURNAL = {Lithuanian Mathematical Journal},
    VOLUME = {49},
      YEAR = {2009},
    NUMBER = {4},
     PAGES = {446--463},
      ISSN = {0363-1672},
   MRCLASS = {60G10 (60F05 60F17 60G52 62M10)},
  MRNUMBER = {2591877},
MRREVIEWER = {Siegfried H\"{o}rmann},
       DOI = {10.1007/s10986-009-9060-x},
       URL = {https://doi.org/10.1007/s10986-009-9060-x},
}

@article {puplinskaite2010,
    AUTHOR = {Puplinskait\.{e}, Donata and Surgailis, Donatas},
     TITLE = {Aggregation of a random-coefficient {${\rm AR}(1)$} process
              with infinite variance and idiosyncratic innovations},
   JOURNAL = {Adv. in Appl. Probab.},
  FJOURNAL = {Advances in Applied Probability},
    VOLUME = {42},
      YEAR = {2010},
    NUMBER = {2},
     PAGES = {509--527},
      ISSN = {0001-8678},
   MRCLASS = {62M10 (60G18 60G52)},
  MRNUMBER = {2675114},
       DOI = {10.1239/aap/1275055240},
       URL = {https://doi.org/10.1239/aap/1275055240},
}

@article{donhauzer2025construction,
author = {Donhauzer, Illia and Leonenko, Nikolai and Olenko, Andriy},
year = {2025},
month = {05},
pages = {},
journal = {ArXiv preprint: 2505.21195},
title = {Construction and limit theorems for sup{CAR} fields},
doi = {10.48550/arXiv.2505.21195}
}

@article {bnleonenko2005,
    AUTHOR = {Barndorff-Nielsen, O. E. and Leonenko, N. N.},
     TITLE = {Spectral properties of superpositions of
              {O}rnstein-{U}hlenbeck type processes},
   JOURNAL = {Methodol. Comput. Appl. Probab.},
  FJOURNAL = {Methodology and Computing in Applied Probability},
    VOLUME = {7},
      YEAR = {2005},
    NUMBER = {3},
     PAGES = {335--352},
      ISSN = {1387-5841},
   MRCLASS = {60E07 (60G10 60G18 62M10)},
  MRNUMBER = {2210585},
MRREVIEWER = {V\'{\i}ctor M. Rivero},
       DOI = {10.1007/s11009-005-4521-0},
       URL = {https://doi.org/10.1007/s11009-005-4521-0},
}

@article {GLT19,
    AUTHOR = {Grahovac, Danijel and Leonenko, Nikolai N. and Taqqu, Murad
              S.},
     TITLE = {Limit theorems, scaling of moments and intermittency for
              integrated finite variance sup{OU} processes},
   JOURNAL = {Stochastic Process. Appl.},
  FJOURNAL = {Stochastic Processes and their Applications},
    VOLUME = {129},
      YEAR = {2019},
    NUMBER = {12},
     PAGES = {5113--5150},
      ISSN = {0304-4149},
   MRCLASS = {60F05 (60G10 60G18 60G22 60G51 60G52)},
  MRNUMBER = {4025701},
MRREVIEWER = {Rita Giuliano Antonini},
       DOI = {10.1016/j.spa.2019.01.010},
       URL = {https://doi.org/10.1016/j.spa.2019.01.010},
}

@article {GLST2019,
    AUTHOR = {Grahovac, Danijel and Leonenko, Nikolai N. and Sikorskii, Alla
              and Taqqu, Murad S.},
     TITLE = {The unusual properties of aggregated superpositions of
              {O}rnstein-{U}hlenbeck type processes},
   JOURNAL = {Bernoulli},
  FJOURNAL = {Bernoulli. Official Journal of the Bernoulli Society for
              Mathematical Statistics and Probability},
    VOLUME = {25},
      YEAR = {2019},
    NUMBER = {3},
     PAGES = {2029--2050},
      ISSN = {1350-7265},
   MRCLASS = {60G51 (60E07 60G10)},
  MRNUMBER = {3961239},
MRREVIEWER = {Nizar Demni},
       DOI = {10.3150/18-BEJ1044},
       URL = {https://doi.org/10.3150/18-BEJ1044},
}

@article {GLT2020multifaceted,
    AUTHOR = {Grahovac, Danijel and Leonenko, Nikolai N. and Taqqu, Murad
              S.},
     TITLE = {The multifaceted behavior of integrated sup{OU} processes: the
              infinite variance case},
   JOURNAL = {J. Theoret. Probab.},
  FJOURNAL = {Journal of Theoretical Probability},
    VOLUME = {33},
      YEAR = {2020},
    NUMBER = {4},
     PAGES = {1801--1831},
      ISSN = {0894-9840},
   MRCLASS = {60F05 (60G10 60G52)},
  MRNUMBER = {4166183},
MRREVIEWER = {Zbigniew Michna},
       DOI = {10.1007/s10959-019-00935-8},
       URL = {https://doi.org/10.1007/s10959-019-00935-8},
}

@article {grahovackevei2025,
    AUTHOR = {Grahovac, Danijel and Kevei, P\'{e}ter},
     TITLE = {Almost sure growth of integrated sup{OU} processes},
   JOURNAL = {Bernoulli},
  FJOURNAL = {Bernoulli. Official Journal of the Bernoulli Society for
              Mathematical Statistics and Probability},
    VOLUME = {31},
      YEAR = {2025},
    NUMBER = {4},
     PAGES = {2597--2623},
      ISSN = {1350-7265},
   MRCLASS = {60E07 (60F15 60G17 60G51 60G57)},
  MRNUMBER = {4931343},
       DOI = {10.3150/24-BEJ1818},
       URL = {https://doi.org/10.3150/24-BEJ1818},
}

@article {Zaffaroni2004,
    AUTHOR = {Zaffaroni, Paolo},
     TITLE = {Contemporaneous aggregation of linear dynamic models in large
              economies},
   JOURNAL = {J. Econometrics},
  FJOURNAL = {Journal of Econometrics},
    VOLUME = {120},
      YEAR = {2004},
    NUMBER = {1},
     PAGES = {75--102},
      ISSN = {0304-4076,1872-6895},
   MRCLASS = {62M10 (62P20)},
  MRNUMBER = {2047781},
       DOI = {10.1016/S0304-4076(03)00207-0},
       URL = {https://doi.org/10.1016/S0304-4076(03)00207-0},
}

\end{document}